\documentclass[10pt]{article}
\usepackage{amsmath,amsfonts,amssymb,euscript,upgreek}

\newtheorem{thm}{Theorem}[section]
\newtheorem{lem}[thm]{Lemma}
\newtheorem{cor}[thm]{Corollary}
\newtheorem{prop}[thm]{Proposition}

        \newcommand{\al}{\alpha}
        \newcommand{\gam}{\gamma}
        \newcommand{\del}{\delta}
        \newcommand{\eps}{\epsilon}
        \newcommand{\veps}{\varepsilon}
        \newcommand{\lam}{\lambda}

        \newcommand{\om}{\omega}
        \newcommand{\tht}{\theta}
        
        \newcommand{\Gam}{{\it \Gamma}}

\newcommand{\RE}{{\mathbb{R}}}
\newcommand{\CO}{{\mathbb{C}}}
\newcommand{\HA}{{\mathbb{H}}}

\newcommand{\HH}{{\EuScript{H}}}
\newcommand{\JJ}{{\EuScript{J}}}
\newcommand{\KK}{{\EuScript{K}}}

\newcommand{\OO}{{\EuScript{O}}}
\newcommand{\VV}{{\EuScript{V}}}
\def\PP{{\mathcal P}}
\def\UU{{\mathcal U}}

\def\A{{\mathsf{A}}}
\def\F{{\mathsf{F}}}
\def\s{{\mathsf{s}}}

\newcommand\nab{\nabla\!}
\newcommand\medwedge{\bigwedge}
\newcommand\sech{\operatorname{sech}}
\newcommand\Pf{\operatorname{Pf}}
\newcommand\End{\operatorname{End}}
\newcommand\Hom{\operatorname{Hom}}
\newcommand\GL{\operatorname{GL}}
\newcommand\SO{\operatorname{SO}}
\newcommand\Spin{\operatorname{Spin}}
\newcommand\Sp{\operatorname{Sp}}
\newcommand\Mp{\operatorname{Mp}}
\newcommand\U{\operatorname{U}}
\newcommand\gl{\mathfrak{gl}}
\newcommand\so{\mathfrak{so}}
\newcommand\rsp{\mathfrak{sp}}
\newcommand\uu{\mathfrak{u}}
\newcommand\kk{\mathfrak{k}}

\newcommand\mm{\mathfrak{m}}
\newcommand\tr{\operatorname{tr}}

\newcommand\sym{\operatorname{Sym}}
\newcommand\id{\mathrm{id}}
\def\eq{&\!\!\!\!=\!\!\!\!&}

\def\qed{$\hfill\Box$}
\newcommand\set[2]{\{#1\,|\,#2\}}
\newcommand\bigset[2]{\Big\{\,#1\,\Big|\,#2\,\Big\}}
\newcommand\mat[4]{\Big(\substack{#1\quad\!#2\vspace{5pt}\\#3\quad\!#4}\Big)}
\def\bra{\langle}
\def\ket{\rangle}
\def\sq{/\!\!/}
\def\pdr{\partial}
\def\ii{\sqrt{-1}\,}
\def\diag{\mathrm{diag}}
\newcommand\textfrac[2]{{\textstyle \frac{#1}{#2}}\,}
\def\hf{\textfrac{1}{2}}
\def\qt{\textfrac{1}{4}}
\def\hfjj{\textfrac{J_0+J_1}{2}}

\begin{document}
\date{}
\title{\vspace{40pt}
Projective flatness in the quantisation of bosons and fermions}
\author{Siye Wu\thanks{Department of Mathematics, University of Hong Kong,
Pokfulam Road, Hong Kong, China (current address)
\newline\hspace*{.175in}
and Department of Mathematics, University of Colorado, Boulder,
CO 80309-0395, USA
\newline\hspace*{.175in}
Email address: {\tt swu@maths.hku.hk}}}

\maketitle

\vspace{-158pt}
\begin{flushright}
{\tt arXiv:1008.5333[math.SG]}\\
HKU-IMR preprint IMR2010:\#11\\
August 2010 (revised: October 2010) 
\end{flushright}

\vspace{100pt}

\begin{abstract}
We compare the quantisation of linear systems of bosons and fermions.
We recall the appearance of projectively flat connection and results on
parallel transport in the quantisation of bosons.
We then discuss pre-quantisation and quantisation of fermions using the
calculus of fermionic variables.
We then define a natural connection on the bundle of Hilbert spaces and
show that it is projectively flat.
This identifies, up to a phase, equivalent spinor representations constructed
by various polarisations.
We introduce the concept of metaplectic correction for fermions and show
that the bundle of corrected Hilbert spaces is naturally flat.
We then show that the parallel transport in the bundle of Hilbert spaces
along a geodesic is the rescaled projection or the Bogoliubov transformation
provided that the geodesic lies within the complement of a cut locus.
Finally, we study the bundle of Hilbert spaces when there is a symmetry.

\medskip\noindent
Keywords: geometric quantisation, projectively flat bundles,
Hermitian symmetric spaces\\
MSC(2000): Primary 53D50; Secondary 32M15, 15A66, 58C50
\end{abstract}

\vspace{10pt}

\section{Introduction}\label{intro}

One of the central questions in geometric quantisation is whether the quantum
Hilbert spaces constructed from different choices of polarisations can be
naturally identified.
Since a quantum state actually corresponds to a ray of vectors, identification
is only required for the projectivisation of the Hilbert spaces.
This amounts to the existence of a natural projectively flat connection on
the bundle of Hilbert spaces over the space of polarisations.
When the symplectic manifold is K\"ahler, it is convenient to consider a
subclass of polarisations that come from complex structures compatible with
the symplectic form.
Given such a complex structure, the quantum Hilbert space is the space of
holomorphic sections of the pre-quantum line bundle.
However, under quite general conditions (satisfied by, for example, the
$2$-sphere), there is no naturally projectively flat connection in the bundle
of quantum Hilbert spaces \cite{GM}.

The next best scenario is that projective flatness holds if we limit the
polarisations to a smaller subset, for example, to those respecting the
symmetry of the system.
For a symplectic vector space polarised by linear complex structures, there
is indeed a natural projectively flat connection in the bundle of Hilbert
spaces \cite{ADPW}.
Moreover, the connection is flat if we include metaplectic correction
\cite{W92,KW}.
Parallel transport in the bundle yields the familiar Fourier and
Segal-Bargmann transforms that are usually used to identify wave functions
in various pictures \cite{KW}.
(The Segal-Bargmann transform can be generalised to relate polarisations on
the cotangent bundles of compact Lie groups \cite{Ha,FMMN}.)
Another example of projective flatness is from quantising the space of flat
connections on a compact orientable surface \cite{Hi,ADPW}; in this case
the complex structures are induced by those on the surface.

Let $(V,\om)$ be a symplectic vector space and $J$, a compatible linear
complex structure on $V$.
The qunatum Hilbert space is a representation of Heisenberg algebra, generated
by tensor powers of $V$ subject to the canonical commutation relation.
The existence of projective flatness is related to the celebrated
Stone-von~Neumann theorem \cite{vN}, which asserts that the irreducible
representation of the Heisenberg algebra is unique up to a unitary equivalence.
Moreover, by Schur's lemma, any two unitary equivalences between two
irreducible representations have to differ only by a phase.
So between the fibres over two linear complex structures in the bundle of
Hilbert spaces, there is a unitary identification which is unique up to a
phase.
This is the hallmark of projectively flat bundles.

The main purpose of this paper is to establish a similar structure of
projective flatness in the quantisation of fermions.
The phase space of a linear fermionic system is a Euclidean space $(V,g)$ and
quantisation means finding an irreducible representation of the Clifford
algebra, which is the fermionic analog of the Heisenberg algebra.
(See \cite{Ko73,Cru} for formal similarities between the two algebras.)
Such a representation is the spinor representation, and just like the bosonic
case, it is unique (when $\dim V$ is even) or nearly unique (when $\dim V$ is
odd) \cite{Ch}.
In the construction of the spinor representation, one needs to choose a
compatible complex structure (see for example \cite{BGV}, \S3.2), which is
the fermionic counterpart of polarisation.
We therefore expect a projectively flat bundle (of spinor representations)
over space of such complex structures.

The rest of the paper is organised as follows.
In \S\ref{bos}, we review the pre-quantisation and quantisation of bosonic
systems whose phase spaces are symplectic vector spaces.
We then recall the natural connection on the bundle of Hilbert spaces and
give a straightforward proof of its projective flatness \cite{ADPW}.
The connection becomes flat after metaplectic correction is included
\cite{W92,KW}.
We present, in a coordinate-free way, the results in \cite{KW} on the parallel
transport in the bundle of Hilbert spaces along geodesics in the base space.
Finally, when there is a group acting symplectically on the vector space, we
decompose the bundle of Hilbert spaces into a direct sum of projectively flat
sub-bundles and identify the invariant part as the bundle from quantising the
symplectic quotient.
\S\ref{ferm} is devoted to the quantisation of fermions when the phase space
is an even-dimensional Euclidean space.
We discuss the pre-quantisation and quantisation of fermions using calculus
of fermionic variables.
We then define a natural connection on the bundle of Hilbert spaces and
show that it is projectively flat.
This identifies, up to a phase, constructions of the spinor representation
under various polarisations.
We introduce the concept of metaplectic correction for fermions and show
that the bundle of corrected Hilbert spaces is naturally flat.
We then show that the parallel transport in the bundle of Hilbert spaces
along a geodesic is the rescaled projection or the Bogoliubov transformation
provided the geodesic lies within the complement of a cut locus.
The decomposition of the bundle of Hilbert spaces when there is a symmetry is
also studied.
In \S\ref{con}, we conclude by highlighting the similarities and differences
in the quantisation of bosons and fermions. 
In Appendix~\ref{A}, we consider the geometry of the spaces of complex
structures compatible to a symplectic or Euclidean structure, which are
classical Hermitian symmetric spaces \cite{Si,H}.
We describe cut locus in the space of polarisations of a fermionic system.
Appendix~B is on the calculus of fermionic variables.
We describe fermionic coherent states and the fermionic analog of Bergman
kernel.
In Appendix~C, we collect some facts on real and quaternionic representations
and on complex structures invariant under a representation.

\section{Quantisation of bosonic systems}\label{bos}

\subsection{Pre-quantisation and quantisation}

Let $(V,\om)$ be a symplectic vector space of dimension $2n$.
A pre-quantum line bundle $\ell$ over $V$ is a line bundle with a connection
whose curvature is $\om/\ii$.
The pre-quantum Hilbert space is $\HH_0=L^2(V,\ell)$, the space of
$L^2$-sections of $\ell$ with respect to the symplectic volume form
$\veps_\om=\om^{\wedge n}/n!$ or $\tilde\veps_\om=\veps_\om/(2\pi)^n$ on $V$.
The covariant derivative $\nabla_x$ along a constant vector field on $V$
parallel to $x\in V$ is a skew-self-adjoint operator on $\HH_0$ and satisfies
the commutation relation $[\nabla_x,\nabla_y]=\om(x,y)/\ii$ for any
$x,y\in V$. 
As $V$ is contractible, $\ell$ is  topologically trivial and is unique up
to an isomorphism.
We can choose a trivialisation of $\ell$ identifying $\HH_0$ with
$L^2(V,\CO)$ such that
\[ \nabla_x=L_x+\hf\ii\iota_x\om,\quad x\in V.  \]
Here $\iota_x\om\in V^*$ is regarded as a linear function on $V$ multiplying
on the sections of $\ell$ or on $L^2(V,\CO)$.
For any $\al\in V^*$, the corresponding pre-quantum operator acting on
$\HH_0$ is
\[ \hat\al=\ii\nabla_{\nu^{-1}(\al)}+\al=\ii L_{\nu^{-1}(\al)}+\hf\al. \]
These operators are self-adjoint on $\HH_0$ and satisfy Heisenberg's canonical
commutation relation $[\hat\al,\hat\beta]=\ii\om^{-1}(\al,\beta)$, where
$\al,\beta\in V^*$.

Consider the space $\JJ_\om$ of compatible complex structures on $(V,\om)$.
(We refer the reader to \S A.1-2 for notations and results on complex
structures.)
For each $J\in\JJ_\om$, the complex subspaces $V^{1,0}_J$, $V^{0,1}_J$ of
$V^\CO$ are Lagrangian with respect to $\om$ and they determine a (linear)
complex polarisation of $(V,\om)$.
The quantum Hilbert space associated to $J$ is
\[  \HH_J=\set{\psi\in\HH_0}{\nabla_x\psi=0,\;\forall x\in V^{0,1}_J}.  \]
So a vector $\psi\in\HH_J$ is a holomorphic $L^2$-section of $\ell$.
The compatibility condition on $J$ guarantees that the space $\HH_J$
is non-empty; in fact, it is infinite dimensional. 
Note that for any $J\in\JJ_\om$, $\HH_J$ is a subspace of $\HH_0$.
Thus we have a bundle of quantum Hilbert spaces $\HH\to\JJ_\om$ whose fibre
over $J\in\JJ_\om$ is $\HH_J$.

The following results are well known.

\begin{prop}\label{bos-holo}
1. Any $\psi\in\HH_J$ is of the form 
\[  \psi=e^{-\frac{1}{4} q_J}\phi   \]
for a unique $J$-holomorphic function $\phi$ on $V$, where 
$q_J\in\om(\cdot,J\cdot)\in\sym^2(V^*)$ is regarded as a quadratic function
on $V$.\\
2. If $\phi$ is a $J$-holomorphic function on $V$, then
$\psi=e^{-\frac{1}{4} q_J}\phi$ is in $\HH_J$ if and only if its norm
\[  \left(\int_Ve^{-\frac{1}{2}q_J}|\phi|^2\,\tilde\veps_\om\right)^{1/2}    \]
weighted by $e^{-\frac{1}{2}q_J}$ is finite, in which case it is equal to
the norm of $\psi\in\HH_0$.\\
3. For any $\al\in V^*$, $\hat\al$ preserves $\HH_J$ and is self-adjoint on
$\HH_J$.
It acts on $\phi$ by
\[  \hat\al\colon\phi\mapsto\ii L_{\nu^{-1}(\al^{0,1})}\phi+\al^{1,0}\phi.  \]
\end{prop}

If $x^i$ ($1\le i\le n$) are the complex coordinates on $V_J^{1,0}$ with
respect to a basis $\{e_i\}_{i=1,\dots,n}$, then the covariant derivative
along $\bar e_{\bar j}$ is 
$\nabla_{\bar j}=\frac{\pdr}{\pdr\bar x^{\bar j}}+\hf q_{i\bar j}x^i$,
where $q_{i\bar j}=q_J(e_i,\bar e_{\bar j})$.
A section $\psi\in\HH_J$ can be identified as a function of the form
\[ \psi(x)=\phi(x)\,\exp[-\hf q_{i\bar j}x^i\bar x^{\bar j}],   \]
where $\phi(x)$ is a holomorphic function in $x=(x^1,\cdots,x^n)\in\CO^n$.

\subsection{Projectively flat connection and metaplectic correction}

The vector bundle $\HH\to\JJ_\om$ of quantum Hilbert spaces is a sub-bundle
of the product bundle  $\JJ_\om\times\HH_0\to\JJ_\om$ of pre-quantum Hilbert
spaces.
The trivial connection on the latter induces a natural connection on $\HH$ by
orthogonal projection.
In \cite{ADPW}, it was shown that the connection on $\HH$ is projectively flat.
For completeness and for comparison with the fermionic case (\S\ref{ferm}.2),
we give a simple derivation of this result.

We first study the effect of the variation $\del J$ on $\HH_J$.
We choose a basis $\{e_i\}_{1\le i\le n}$ of $V_J^{1,0}$.
Suppose $\psi\in\HH_J$, i.e., $\nabla_{\bar k}\psi=0$ ($1\le k\le n$).
As $J$ changes to $J+\del J$, the infinitesimal parallel transport
$\psi+\del\psi\in\HH_{J+\del J}$ of $\psi\in\HH_J$ is the orthogonal
projection of $\psi$ on $\HH_{J+\del J}$.
Thus we have two conditions: $\del\psi\perp\HH_J$, i.e.,
$\del\psi\perp\ker\nab_{\bar i}$ for $1\le i\le n$, and
$\psi+\del\psi\in\HH_{J+\del J}$, or
 \[  \nabla_{\bar e_{\bar i}+\del\bar e_{\bar i}}(\psi+\del\psi)=0,  \]
where $\del\bar e_{\bar i}=(\del\bar P)_{\bar i}^{\;\,j}e_j$ (see \S\ref{A}.2).
This implies, to the first order, that $\del\psi$ satisfies the equation
\[ \nab_{\bar i}\,(\del\psi)=-(\del\bar P)_{\bar i}^{\;\,j}\nab_j\psi
   =(\del P)_{\bar i}^{\;\,j}\nab_j\psi.   \]
We claim that
\[  \del\psi=\hf\ii\nab_i(\del P)^{ij}\nab_j\psi
            =\hf\ii(\del P)^{ij}\nab_i\nab_j\psi  \]
is the (unique) solution satisfying the above conditions.
(The second equality holds because $(\del P)^{ij}$ is a constant tensor 
on $V$.)
First, this $\del\psi$ is orthogonal to $\HH_J$ as $\nab_{\bar i}$ is the
formal adjoint of $\nab_i$.
Second, as $\nab_{\bar i}\psi=0$, $(\del P)^{ij}=(\del P)^{ji}$ and
$[\nab_{\bar i},\nab_j]=\om_{\bar ij}/\ii$, we get
\begin{eqnarray*}
\nab_{\bar k}\,(\del\psi)\eq\hf\ii(\del P)^{ij}[\nab_{\bar k},\nab_i\nab_j]\psi
=\hf\ii(\del P)^{ij}([\nab_{\bar k},\nab_i]\nab_j
                                  +\nab_i[\nab_{\bar k},\nab_j])\psi \\
\eq\om_{\bar ki}(\del P)^{ij}\nab_j\psi=(\del P)_{\bar k}^{\;\,j}\nab_j\psi.
\end{eqnarray*}
The uniqueness is clear from the geometric interpretation.

The connection $1$-form $\A^\HH$ on $\HH$ satisfies $(\del+\A^\HH)\psi=0$.
Therefore $\A^\HH=-\frac{\ii}{2}(\del P)^{ij}\nab_i\nab_j$;
it is an operator-valued $1$-form on $\JJ_\om$.
We then calculate
\begin{eqnarray*}
\del\,\A^\HH\eq\textfrac{\ii}{2}(\del P)^{ij}\wedge\nabla_{\del e_i}\nabla_j
    =\textfrac{\ii}{2}(\del P)^{ij}\wedge(\del P)_i^{\;\bar k}
     \nabla_{\bar k}\nabla_j \\
    \eq\hf\om_{\bar kj}(\del P)^{ji}\wedge(\del P)_i^{\;\bar k}
    =\hf\tr(P\,\del P\wedge\del P\,P),
\end{eqnarray*}
ignoring the terms that vanish on $\HH_J$, and
\[  \A^\HH\wedge\A^\HH
 =-\qt\nabla_i\nabla_j\nabla_k\nabla_l\,(\del P)^{ij}\wedge(\del P)^{kl}=0. \]
Therefore the curvature of the connection on $\HH$ is 
\[  \F^\HH=\hf\tr(P\,\del P\wedge\del P\,P)\,\id_\HH
          =\upsigma_\om/2\ii\,\id_\HH.  \]
Since it is a $2$-form on $\JJ$ times the identity operator on the fibre,
the connection on $\HH$ is indeed projectively flat \cite{ADPW}.
(In \S2.1 of \cite{W93}, it was shown directly, without the orthogonal
projection from $\HH_0$, that the formula for $\A^\HH$ defines a connection
on $\HH$ which is projectively flat.)

We now incorporate metaplectic correction. 
Denote the restriction of $\KK=(\det\!\VV)^*\to\JJ$ to $\JJ_\om$
(see \S\ref{A}.1) by the same notation $\KK$.
Since $\JJ_\om$ is contractible, there is a unique bundle
$\sqrt\KK\to\JJ_\om$ such that $(\sqrt\KK)^{\otimes2}=\KK$.
The curvature of the natural connection on $\sqrt\KK$ is
(see \cite{KW} or \S\ref{A}.1)
\[    \F^{\sqrt\KK}=\hf\F^\KK=-\hf\tr(P\,\del P\wedge\del P\,P)
                   =-\upsigma_\om/2\ii.  \]
We consider the bundle $\hat\HH=\HH\otimes\sqrt\KK$.
Its fibre $\hat\HH_J=\HH_J\otimes\sqrt{\KK_J}$ over $J\in\JJ_\om$ is the
metaplectically corrected quantum Hilbert space with the polarisation $J$.
Since the curvatures of $\HH$ and $\sqrt\KK$ cancel, the bundle
$\hat\HH\to\JJ_\om$ is canonically flat \cite{W92,KW}.
The flatness of the bundle indicates that for the symplectic linear space,
quantisation is independent of the choice of polarisations. 
We summarise the results in the following

\begin{thm} {\rm (\cite{ADPW,W92,KW})}
Consider the quantisation of a bosonic system whose phase space is a
finite dimensional symplectic vector space $(V,\om)$.\\
1. The bundle of quantum Hilbert spaces $\HH\to\JJ_\om$ is projectively 
flat, with curvature $\upsigma_\om/2\ii$.\\
2. The bundle of quantum Hilbert spaces with metaplectic correction
$\hat\HH\to\JJ_\om$ is flat.
\end{thm}

\subsection{Parallel transport along geodesics and
the Bogoliubov transformations}

We recall various results in \cite{KW}.
Let $(V,\om)$ be a symplectic vector space.
The space $(\JJ_\om,\upeta_\om)$ of compatible complex structures is
non-positively curved and there is a unique geodesic connecting any two points.
Let $J_0,J_1\in\JJ_\om$ define two complex polarisations.
We want to study the parallel transport $\UU_{J_1J_0}^\HH$ and
$\UU_{J_1J_0}^{\hat\HH}$ in the bundles $\HH$ and $\hat\HH$, respectively,
along the geodesic from $J_0$ to $J_1$.
A related notion is the orthogonal projection $\PP_{J_1J_0}$ from
$\HH_{J_0}$ to $\HH_{J_1}$ in $\HH_0$.

\begin{thm}\label{bos-bog} {\rm (\cite{KW})}
Let $J_0,J_1\in\JJ_\om$ and let $\gam=\{J_t\}_{0\le t\le1}$ be the (unique)
geodesic from $J_0$ to $J_1$, $t$ being proportional to the arc-length
parameter.
Then\\
1. the parallel transport in $\HH$ along $\gam$ is 
$\UU_{J_1J_0}^\HH=\left(\det\hfjj\right)^{1/4}\,\PP_{J_1J_0}$;\\
2. the parallel transport in $\sqrt\KK$ along $\gam$ is
$\UU_{J_1J_0}^{\sqrt\KK}=\sqrt{\;\;\big(\det J_{1/2}/\ii\big|_{V_{J_0}^{1,0}}
\big)^{-1\hspace{-106pt}T\hspace{100.25pt}}}$, and 
\[ \bra\,\UU_{J_1J_0}^{\sqrt\KK}\sqrt{\mu'_0},\sqrt{\mu_0}\,\ket
=\left(\det\hfjj\right)^{-1/4}\,\bra\sqrt{\mu'_0},\sqrt{\mu_0}\,\ket  \]
for any $\mu_0,\mu'_0\in\medwedge^n(V_{J_0}^{1,0})^*$;\\
3. the parallel transport in $\hat\HH$ along $\gam$, which is
$\UU_{J_1J_0}^{\hat\HH}=\UU_{J_1J_0}^\HH\otimes\UU_{J_1J_0}^{\sqrt\KK}$,
corresponds to the pairing between $\hat\HH_{J_0}$ and $\hat\HH_{J_1}$ given by
\[ \bra\psi_1\otimes\sqrt{\mu_1},\psi_0\otimes\sqrt{\mu_0}\,\ket
   =\bra\psi_1,\psi_0\ket\bra\sqrt{\mu_1},\sqrt{\mu_0}\,\ket, \qquad
    \psi_l\in\HH_{J_l},\;\mu_l\in{\textstyle \medwedge^n}
    (V_{J_l}^{1,0})^*\;\,(l=0,1). \]
\end{thm}

\noindent{\em Proof:}
Part~1 is Theorem~3.4 of \cite{KW}.
Part~2 follows from Theorem~3.3.2 and formula~(3.9) of \cite{KW}, except
the parallel transport itself is expressed more intrinsically using
Proposition~\ref{J/2}.
Part~3 is Corollary~3.7 of \cite{KW}.
\qed\vspace{10pt}

We note that $\hfjj$ is always invertible for $J_0,J_1\in\JJ_\om$.
The factor $\left(\det\hfjj\right)^{1/4}$ appeared in \cite{Da,W81} and
was used to rescale the projection $\PP_{J_1J_0}$ to a unitary operator
called the Bogoliubov transformation \cite{W81,W92}.
Therefore Theorem~\ref{bos-bog}.1 shows that the parallel transport
$\UU_{J_1J_0}^\HH$ along the geodesic coincides with the Bogoliubov
transformation.
(The induced parallel transport on the creation and annihilation operators
gives the more traditional version of Bogoliubov transformations.) 
The parallel transport can also be expressed using the Bergman kernel
(Proposition~3.6 of \cite{KW}):

\begin{cor}\label{bos-berg} {\rm (\cite{KW})}
Let $\psi=\phi\,e^{-\frac{1}{4}q_{J_0}}\in\HH_{J_0}$.
Then for $x\in V$,
\[  (\UU_{J_1J_0}^\HH\psi)(x)=\left(\det\hfjj\right)^{1/4}
  \,e^{-\frac{1}{4}q_{J_1}(x)}\int_V\exp[\ii\om(x_{J_1}^{1,0},y)
  -\qt q_{J_1}(y)-\qt q_{J_0}(y)]\,\phi(y)\;\tilde\veps_\om(y).     \]
\end{cor}

Of particular interest is the parallel transport of a coherent state
\[     c_J^\al(x)=\exp[q_J(\bar\al,x)-\qt q_J(x)]
       =\exp[\ii\om(\bar\al,x_J^{1,0})-\qt q_J(x)],\quad x\in V,  \]
where $J\in\JJ_\om$ and $\al\in V_J^{1,0}$ is a parameter.
We recall some results from \cite{KW}, but in a coordinate-free way.

\begin{thm}\label{bos-coh} {\rm (\cite{KW})}
Under the assumptions of Theorem~\ref{bos-bog},\\
1. the parallel transport along $\gam$ of the coherent state
$c_{J_0}^\al\in\HH_{J_0}$ is, for $x\in V$,
\[  (\UU^\HH_{J_1J_0}c_{J_0}^\al)(x)
    =\left(\det\hfjj\right)^{-1/4}e^{-\frac{1}{4}q_{J_1}(x)}
    \exp\big[\hf\om\big(x_{J_1}^{1,0}-\bar\al,
    \,\left(\hfjj\right)^{\!\!-1}\!(x_{J_1}^{1,0}-\bar\al)\,\big)\big];   \]
2. the parallel transport along $\gam$ of any state
$\psi=e^{-\frac{1}{4}q_{J_0}}\phi\in\HH_{J_0}$ is, for $x\in V$,
\[ (\UU_{J_1J_0}^\HH\psi)(x)=\left(\det\hfjj\right)^{\!-1/4}
   e^{-\frac{1}{4}q_{J_1}(x)}\!\int_V\exp\big[\hf\om\big(
   x_{J_1}^{1,0}-y_{J_0}^{0,1},\,\left(\hfjj\right)^{\!\!-1}\!
   (x_{J_1}^{1,0}-y_{J_0}^{0,1})\,\big)-\hf q_{J_0}(y)\big]
   \,\phi(y)\;\tilde\veps_\om(y).\]
\end{thm}

\noindent{\em Proof:}
Part~1 is Theorems~3.3.1 of \cite{KW}, where it was proved by solving the
equation of parallel transport.
As remarked in \cite{KW} (after Corollary~3.7), the result also follows from
the Bergman kernel by Theorem~\ref{bos-bog}.1 or Corollary~\ref{bos-berg}.
Since the latter approach will be adapted in the proof of
Theorem~\ref{ferm-coh}.1 for fermions, we include the details here
for comparison.
Indeed, for any $x\in V$,
\begin{eqnarray*}
(\PP_{J_1J_0}c^\al_{J_0})(x)\eq
 e^{-\frac{1}{4}q_{J_1}(x)}\int_V\exp\big[\ii\om(y,x_{J_1}^{1,0})
 -\hf\om\big(y,\hfjj y\big)\big]\,e^{\ii\om(\bar\al,y)}\;\tilde\veps_\om(y) \\
\eq e^{-\frac{1}{4}q_{J_1}(x)}
 \exp\big[\hf\om\big(x_{J_1}^{1,0}-\bar\al,\,
 \big(\hfjj\big)^{\!\!-1}\!(x_{J_1}^{1,0}-\bar\al)\big)\big]
 \int_Ve^{-\hf\om\big(y',\hfjj y'\big)}\;\tilde\eps(y')   \\
\eq\left(\det\hfjj\right)^{-1/2}e^{-\frac{1}{4}q_{J_1}(x)}
 \exp\big[\hf\om\big(x_{J_1}^{1,0}-\bar\al,\,
 \big(\hfjj\big)^{\!\!-1}\!(x_{J_1}^{1,0}-\bar\al)\big)\big],
\end{eqnarray*}
where the change of variable is $y'=y-\ii(\hfjj)^{-1}(x^{1,0}_{J_1}-\bar\al)$.
Here we used the Gaussian integral
\[  \int_V e^{-\frac{1}{2}\om(x,Ax)}\,\tilde\veps_\om(x)=(\det A)^{-1/2}   \]
for any $A\in\End(V)$ such that $\om(\cdot,A\cdot)$ is a symmetric,
positive-definite bilinear form.
(This implies $\det A>0$.)
Part~2 is Theorems~3.8 of \cite{KW}.
\qed\vspace{10pt}

In particular,  when $J_1=J_0$, the above reduces to the identity
transformation.
When $n=1$, we use the parametrisation in \S A.2.
We note that $\al\in V_{J_0}^{1,0}$ and $x\in V\cong V_{J_1}^{1,0}$ can be
identified with complex numbers.
If the geodesic $\{J_t\}_{0\le t\le 1}$ from $J_0$ to $J_1$ is given by
$z(t)=\tanh bt$, then Theorem~\ref{bos-coh} gives
\[ (\UU^\HH_{J_1J_0}c_{J_0}^\al)(x)=\sqrt{\sech b}\,
\exp\big[\,\bar\al x\tanh b+\hf(\bar\al^2-x^2)\sech b-\qt|x|^2\,\big].   \]
This is an important case ($n=1$) of Proposition~3.2 in \cite{KW}.

\subsection{Bosonic systems with symmetries}

Let $(V,\om)$ be a symplectic vector space of dimension $2n$.
The action of $\Sp(V,\om)$ on $\JJ_\om$ can be lifted to $\HH$, preserving
the projectively flat connection.
It can be lifted to an action of the metaplectic group $\Mp(V,\om)$, which
is a double cover of $\Sp(V,\om)$, on $\sqrt\KK$ and hence on $\hat\HH$.
The lifted action preserves the flat connection.

Let $K$ be a compact Lie group with Lie algebra $\kk$.
Suppose there is a representation of $K$ on $V$ preserving $\om$.
Then it also acts on the bundles $\HH$ and $\hat\HH$ preserving the
connections.
Over the fixed-point set $(\JJ_\om)^K$, the group $K$ acts on the fibres
of $\HH$.
Each fibre splits orthogonally into a direct sum of subspaces of various
representation types.
Since $K$ preserves the connection, the restriction of the bundle $\HH$ to
$(\JJ_\om)^K$, together with the projectively flat connection, splits into
sub-bundles with fibre-wise $K$-actions.
The sub-bundle $\HH^K$ on which $K$ acts trivially is related to the
quantisation of the symplectic quotient.

The action of $K$ on $V$ is Hamiltonian with a moment map
$\mu_K\colon V\to\kk^*$ given by
\[ \bra\mu_K(x),A\ket=\hf\om(x,Ax),\quad x\in V,\;A\in\kk.   \]
The symplectic quotient $V\sq K=\mu_K^{-1}(0)/K$ is a stratified symplectic
space \cite{SL}.
Choosing $J\in(\JJ_\om)^K$, the action of $K$ extends to that of $K^\CO$.
Let $\pi\colon V\to V/K^\CO$ be the quotient map.
With the above moment map, every point in $V$ is semi-stable, i.e.,
$V^{\mathrm{ss}}=V$ (see Example~2.3 of \cite{Sj}).
The quotient $V/K^\CO=V\sq K$ is also a stratified analytic space; a function
$f$ on an open set $U\subset V/K^\CO$ is analytic if $\pi^*f$ is so on
$\pi^{-1}(U)$.
On the singular space $V\sq K$, this analytic structure replaces the notion
of polarisation.
The sheaf of invariant sections $\pi^K_*\ell$ on $V/K^\CO$ defined by
$\OO(\pi^K_*\ell)(U)=\Gam(\pi^{-1}(U),\OO(\ell))^K$ plays the role of
a pre-quantum line bundle.
We have (cf.~Proposition~2.14 and Theorem~2.18 of \cite{Sj})
\[  \Gam(V\sq K,\OO(\pi^K_*\ell))\cong\Gam(V,\OO(\ell))^K.    \]
Restricting to the $L^2$-subspaces, we can identify $(\HH_J)^K$ with the
quantum Hilbert space arising from the quantisation of $V\sq K$ with a
complex structure induced from $J$.

We have a projectively flat bundle $\HH^K\to(\JJ_\om)^K$ whose fibres are
quantum Hilbert spaces of $V\sq K$ with complex structures from $(\JJ_\om)^K$.
The connection is unitary if the inner product in the fibres $(\HH_J)^K$ is the
restriction of that in $\HH_J$.
This is the case, for example, in the quantisation of Chern-Simons gauge
theory \cite{ADPW}.
The inner product on $(\HH_J)^K$ from quantisation of $V\sq K$ is usually
different.
In \cite{HK}, it was shown that for a compact symplectic manifold with
a Hamiltonian group action and when metaplectic correction is included,
the two inner products agree in the semi-classical limit.

Unless the moment map $\mu_K$ is proper, the symplectic quotient $V\sq K$
is non-compact and the quantum Hilbert space $(\HH_J)^K$ is expected to
be infinite dimensional.
When $\mu_K$ is proper however, the base space of the bundle
$\HH^K\to(\JJ_\om)^K$ is a point.

\begin{prop}\label{proper}
If $\mu_K\colon V\to\kk^*$ is proper, then $(\JJ_\om)^K=\{J_0\}$.
\end{prop}

\noindent{\em Proof:}
If $(\JJ_\om)^K\ne\{J_0\}$, then by Proposition~\ref{K-b-f}.1, there is a
non-zero $K$-invariant complex subspace $(V',J_0)$ of $(V,J_0)$ and a
$K$-invariant real structure $R$ on $V'$ such that $\om(Rx,Ry)=-\om(x,y)$
for all $x,y\in V'$.
For any $A\in\kk$, $x\in V'$, we have
\[  \bra\mu_K(Rx),A\ket=\hf\om(Rx,ARx)=\hf\om(Rx,RAx)
    =-\hf\om(x,Ax)=-\bra\mu_K(x),A\ket.                            \]
Therefore $\mu_K=0$ on $V'_0=(V')^R$ and hence $\mu_K$ is not proper.
\qed\vspace{5pt}

\section{Quantisation of fermionic systems}\label{ferm}

\subsection{Pre-quantisation and quantisation}

We consider pre-quantisation \cite{Ko77} and quantisation \cite{W81} of linear
fermionic systems.
The phase space is given by a finite-dimensional real vector space $V$ equipped
with a Euclidean inner product $g$.
More precisely, it is a fermionic copy $\Pi V$ of $V$ (see \S\ref{B}.1).
The pre-quantum line bundle does not exists in the usual sense, but its
``sections'' and the operators acting on them do.
Motivated by the bosonic case (\S\ref{bos}.1), we take the pre-quantum
Hilbert space $\HH_0$ of fermions as $\medwedge^\bullet(V^\CO)^*$,
the space of ``functions" on $\Pi V$.
On $\HH_0$, there is an Hermitian form given by the Berezin integral
(see \S\ref{B}.1 for definition and notations)
\[  \bra\psi,\psi'\ket_0=\int_{\Pi V}\!\bar\psi\wedge\star_0\,\psi'\;\eps_g,
    \quad\psi,\psi'\in\HH_0,  \]
where $\bar\psi$ is the standard complex conjugation of $\psi$, $\star_0$
is the Hodge star defined by the metric $\hf g$.
The covariant derivatives take the form
\[  \nabla_x=\iota_x-\hf\nu(x)\wedge\,\cdot\,,\quad x\in V  \]
and satisfy the relation
\[   \{\nabla_x,\nabla_y\}=-g(x,y),\quad x,y\in V.   \]
So the ``curvature'' is a symmetric bilinear form; the minus sign is enforced
by the requirement, as in the bosonic case, that the covariant derivatives are
skew-self-adjoint operators on $\HH_0$.

A linear functional $\al\in V^*$ is a ``classical observable" that can be
pre-quantised, giving rise to a self-adjoint operator
\[  \hat\al=\nabla_{\nu^{-1}(\al)}+\al\wedge\,\cdot\,
=\iota_{\nu^{-1}(\al)}+\hf\al\wedge\,\cdot\  \]
on $\HH_0$.
These operators satisfy the canonical anti-commutation relation or 
the Clifford algebra relation
\[   \{\hat\al,\hat\beta\}=g^{-1}(\al,\beta),\quad\al,\beta\in V^*,  \]
making $\HH_0$ a (reducible) Clifford module.

We now assume that $V$ is even dimensional; let $\dim V=2n$.
Recall from \S A.2 the space $\JJ_g$ of complex structures on $V$ compatible
with the metric $g$ and the orientation.
Each $J\in\JJ_g$ defines a polarisation, a maximally isotropic complex
subspace $V_J^{1,0}$ of $V^\CO$.
The quantum Hilbert space (with the choice of polarisation $J$) is
\[  \HH_J=\set{\psi\in\HH_0}{\nabla_x\psi=0,\;\forall x\in V^{0,1}_J}.  \]
We have a bundle of quantum Hilbert spaces $\HH\to\JJ_g$ whose fibre over
$J\in\JJ_g$ is $\HH_J$.

On $\HH_0$, there is an involution $\psi\mapsto\psi^*$ defined as the unique
linear extension of the operation
$(\al_1\wedge\cdots\wedge\al_k)^*=\bar\al_k\wedge\cdots\wedge\bar\al_1$,
where $\al_1,\dots,\al_k\in(V^\CO)^*$.

\begin{prop}\label{ferm-holo}
1. Any $\psi\in\HH_J$ is of the form 
\[  \psi=e^{\frac{\ii}{2}\varpi_J}\wedge\phi   \]
for a unique $\phi\in\medwedge^\bullet(V^{1,0}_J)^*$, where
$\varpi_J=g(J\cdot,\cdot)\in\medwedge^2V^*$.
Consequently, $\dim_\CO\HH_J=2^n$.\\
2. Suppose $\psi,\psi'\in\HH_J$ correspond to
$\phi,\phi'\in\medwedge^\bullet(V^{1,0}_J)^*$, respectively, then
\[  \bra\psi,\psi'\ket_0=\int_{\Pi V}\!\bar\phi\wedge\star\,\phi'\;\eps_g
    =\int_{\Pi V}\phi^*\wedge\phi'\wedge e^{\ii\varpi_J}\,\tilde\eps_g,  \]
where $\star$ is the Hodge star defined by $g$ and
$\tilde\eps_g=\ii^n\,\eps_g$.\\
3. For any $\al\in V^*$, $\hat\al$ preserves $\HH_J$ and remains self-adjoint
on $\HH_J$.
It acts on $\phi\in\medwedge^\bullet(V^{1,0}_J)^*$ by
\[  \hat\al\colon\phi\mapsto
    \iota_{\nu^{-1}(\al^{0,1})}\phi+\al^{1,0}\wedge\phi.  \]
\end{prop}

\noindent{\em Proof:}
1. Write $\psi=e^{\frac{\ii}{2}\varpi_J}\wedge\phi$ for some (unique)
$\phi\in\HH_0$.
Then for any $x\in V$, we have
\[  \nabla_x\psi=e^{\frac{\ii}{2}\varpi_J}\wedge
                 (\iota_x\phi-\nu(x^{1,0})\wedge\phi).  \]
Therefore $\psi\in\HH_J$ if and only if $\iota_x\phi=0$ for all
$x\in V^{0,1}_J$.
This implies $\phi\in\medwedge^\bullet(V^{1,0}_J)^*$.\\
2. We choose a basis $\{e_i\}$ of $V^{1,0}_J$ such that
$g(e_i,\bar e_j)=\del_{ij}$ and $\eps_g=e_1\wedge\cdots\wedge e_n$.
Assume, without loss of generality, that
$\phi=\phi'=e_1^*\wedge\cdots\wedge e_k^*$ ($1\le k\le n$).
Then 
\[ \psi=\phi\wedge\sum_{r=0}^{n-k}\;2^{-r}\!\!
        \sum_{k+1\le i_1<\cdots<i_r\le n}\!\!\bar e_{i_1}^*\wedge
        e_{i_1}^*\wedge\cdots\wedge\bar e_{i_r}^*\wedge e_{i_r}^*. \]
Since $\star_0=2^{p-n}\star$ on $\medwedge^pV^*$, we have
\[ \star_0\,\psi=\frac{(-1)^{\frac{k(k-1)}{2}}}{\ii^n}\,
 \phi\wedge\sum_{r=0}^{n-k}(-1)^r\,2^{r+k-n}\!\!\!\!\!\!
 \sum_{k+1\le j_1<\cdots<j_{n-r}\le n}\!\!\bar e_{j_1}^*\wedge e_{j_1}^*\wedge
 \cdots\wedge\bar e_{j_{n-k-r}}^*\wedge e_{j_{n-k-r}}^*.   \]
So
\[  \bra\psi,\psi'\ket_0=\int_{\Pi V}\!\bar\psi\wedge\star_0\,\psi\;\eps_g
    =\sum_{r=0}^{n-k}2^{k-n}{n-k\choose r}=1. \]
The two integrals in the equality are clearly $1$
(cf.~proof of Proposition~\ref{B}.2).\\
3. This follows from the identity
\[  \hat\al(e^{\frac{\ii}{2}\varpi_J}\wedge\phi)=e^{\frac{\ii}{2}\varpi_J}
    \wedge(\iota_{\nu^{-1}}(\al)\phi+\al^{1,0}\wedge\phi),   \]
which yields the result when $\phi\in\medwedge^\bullet(V^{1,0}_J)^*$.
\qed\vspace{5pt}

We note that the space $\medwedge^\bullet(V^{1,0}_J)^*$ with the action
of $\al\in V^*$ in Proposition~\ref{ferm-holo}.3 is the standard
construction of the irreducible Clifford module of spinors.
Here it arises naturally in the quantisation of fermionic systems.
The factor $e^{\frac{\ii}{2}\varpi_J}$ is the fermionic analogue of the
Gaussian in Proposition~\ref{bos-holo}.1.
It is crucial in achieving projective flatness for the bundle $\HH\to\JJ_g$
(\S\ref{ferm}.2), as the bundle $\medwedge^\bullet\VV^*\to\JJ_g$ without the
fermionic Gaussian factor is not projectively flat.

The results in this section can be explained using ``fermionic coordinates".
We refer the reader to \S\ref{B}.2 where this is done.

\subsection{Projectively flat connection and metaplectic correction}

We study the geometry of the bundle $\HH\to\JJ_g$ of Hilbert spaces of the 
fermionic system $(V,g)$.
Following \S\ref{bos}.2, we define a connection on $\HH$ by orthogonal
projection of the trivial connection on the product bundle
$\JJ_g\times\HH_0\to\JJ_g$.
We now show that this connection is also projectively flat.

Along a variation $\del J$ of $J\in\JJ_g$, $\psi\in\HH_J$ changes to
$\psi+\del\psi\in\HH_{J+\del J}$ by parallel transport.
Since $\psi+\del\psi$ is the (infinitesimal) orthogonal projection of $\psi$
to $\HH_{J+\del J}$, we have, as in \S\ref{bos}.2, $\del\psi\perp\HH_J$ and
\[ \nab_{\bar i}(\del\psi)=-(\del\bar P)_{\bar i}^{\;\,j}\nab_j\psi
   =(\del P)_{\bar i}^{\;\,j}\nab_j\psi.   \]
The (unique) solution that satisfies the above two conditions is
\[  \del\psi=-\hf\nab_i(\del P)^{ij}\nab_j\psi
            =-\hf(\del P)^{ij}\nab_i\nab_j\psi.  \]
First, this $\del\psi$ is orthogonal to $\HH_J$ as $\nab_{\bar i}$ is the
formal adjoint of $\nab_i$.
Second, as $\nab_{\bar i}\psi=0$, $(\del P)^{ij}=-(\del P)^{ji}$ and
$\{\nab_{\bar i},\nab_j\}=-g_{\bar ij}$, we get
\begin{eqnarray*}
\nab_{\bar k}(\del\psi)\eq-\hf(\del P)^{ij}[\nab_{\bar k},\nab_i\nab_j]\psi
 =-\hf(\del P)^{ij}
   (\{\nab_{\bar k},\nab_i\}\nab_j-\nab_i\{\nab_{\bar k},\nab_j\})\psi     \\
  \eq g_{\bar ki}(\del P)^{ij}\nab_j\psi=(\del P)_{\bar k}^{\,\;j}\nab_j\psi.
\end{eqnarray*}
The connection $\A^\HH$ on $\HH$ is determined by $(\del+\A^\HH)\psi=0$
and thus $\A^\HH=\hf(\del P)^{ij}\nab_i\nab_j$.
Following the calculations in \S\ref{bos}.2, we get
\begin{eqnarray*}
\del\,\A^\HH\eq-\hf(\del P)^{ij}\wedge\nabla_{\del e_i}\nabla_j
   =-\hf(\del P)^{ij}\wedge(\del P)_i^{\;\bar k}\nabla_{\bar k}\nabla_j   \\
                  \eq-\hf g_{\bar kj}(\del P)^{ji}\wedge(\del P)_i^{\;\bar k}
   =-\hf\tr(P\,\del P\wedge\del P\,P)
\end{eqnarray*}
and
\[ \A^\HH\wedge\A^\HH=\qt\nabla_i\nabla_j\nabla_k\nabla_l\,
            (\del P)^{ij}\wedge(\del P)^{kl}=0.                            \]
Therefore the curvature of the connection $\A^\HH$ is 
\[  \F^\HH=-\hf\tr(P\,\del P\wedge\del P\,P)\,\id_\HH
          =\upsigma_g/2\ii\,\id_\HH.                       \]
Since it is a $2$-form on $\JJ_g$ times the identity operator on the fibre,
the connection is projectively flat.

We propose a metaplectic correction for fermions. 
Recall the line bundle $\KK^{-1}=\det\VV\to\JJ_g$ whose fibre over $J$ is
$\KK^{-1}_J=\medwedge^nV^{1,0}_J$.
We claim that $c_1(\KK)$ is even.
This can be seen from the holonomy of the bundle $\HH$ with curvature
$\F^\HH=-\frac{1}{2}\F^\KK\;\id_\HH$.
Since $\JJ_g$ is simply connected, there is a unique line bundle
$\sqrt{\KK^{-1}}\to\JJ_g$ such that $(\sqrt{\KK^{-1}})^{\otimes2}=\KK^{-1}$.
The bundle $\sqrt{\KK^{-1}}$ has a connection (\S\ref{A}.1) whose curvature is
\[    \F^{\sqrt{\KK^{-1}}}=-\hf\F^\KK=\hf\tr(P\,\del P\wedge\del P\,P)
                          =-\upsigma_g/2\ii.            \]
For any $J,J'\in\JJ_g$, there is a pairing between 
$\KK_J^{-1}=\medwedge^nV_J^{1,0}$ and $\KK_{J'}^{-1}=\medwedge^nV_{J'}^{1,0}$.
For any $\mu\in\KK_J^{-1}$ and $\mu'\in\KK_{J'}^{-1}$, $\bra\mu',\mu\ket$ is
the ratio of $\bar\mu'\wedge\mu$ and $\tilde\eps_g$.
Since $\bra\mu,\mu\ket>0$ if $\mu\ne0$, there is an inner product on
$\sqrt{{\KK_J}^{-1}}$ defined by
$\bra\sqrt\mu,\sqrt\mu\ket=\sqrt{\bra\mu,\mu\ket}$.

We consider the bundle $\hat\HH=\HH\otimes\sqrt{\KK^{-1}}$.
The fibre $\hat\HH_J=\HH_J\otimes\sqrt{\KK_J^{-1}}$ over $J\in\JJ_g$ is
called the metaplectically corrected quantum Hilbert space in polarisation $J$.
Since the curvatures of $\HH$ and $\sqrt{\KK^{-1}}$ cancel, the bundle
$\hat\HH\to\JJ$ is canonically flat.
The flatness of the bundle indicates that for the fermionic system whose
phase space is a linear space, quantisation is independent of the choice
of polarisations.
We note here that in contrast to the bosonic case, the metaplectic correction
is obtained by tensoring $\HH_J$ with $\sqrt{\KK_J^{-1}}$ instead of
$\sqrt{\KK_J}$.
This is clearly related to the opposite way a fermionic measure transforms
under coordinate changes.
We recall from \S\ref{A}.2 that the pseudo-K\"ahler form $\upsigma$ restricts
to $\upsigma_\om$ on $\JJ_\om$ but to $-\upsigma_g$ on $\JJ_g$.
Consequently, for both bosonic and fermionic systems, the line bundle of
half-forms has a negative first Chern form on the space of polarisations.

We summarise the results in the following

\begin{thm}
Consider the quantisation of a fermionic system whose phase space is given by
a finite dimensional Euclidean vector space $(V,g)$.\\
1. The bundle of quantum Hilbert spaces $\HH\to\JJ_g$ is projectively flat,
with curvature $\upsigma_g/2\ii$.\\
2. The bundle of quantum Hilbert spaces with metaplectic correction
$\hat\HH\to\JJ_g$ is flat.
\end{thm}

\subsection{Parallel transport along geodesics and
the Bogoliubov transformations}

We investigate the parallel transport in the bundles $\HH$ and $\hat\HH$ along
geodesics in $\JJ_g$.
Unlike $\JJ_\om$, which is contractible and non-positively curved, the space
$\JJ_g$ is compact and non-negatively curved.
The geodesics through two conjugate points in $\JJ_g$ are not unique.
Nevertheless, we show that if $J_0,J_1\in\JJ_g$ are not in the cut loci (see
\S\ref{A}.3) of each other, then the parallel transport $\UU_{J_1J_0}^\HH$
along the (unique) length-minimising geodesic from $J_0$ to $J_1$ is the 
rescaled orthogonal projection $\PP_{J_1J_0}$ from $\HH_{J_0}$ to $\HH_{J_1}$
in $\HH_0$.
The latter was called the Bogoliubov transformation of fermionic systems
\cite{W81}.
The inner product in $\sqrt{\KK_{J_0}^{-1}}$ extends to a pairing between
$\sqrt{\KK_{J_0}^{-1}}$ and $\sqrt{\KK_{J_2}^{-1}}$ along the geodesic,
which is non-degenerate as long as $J_1$ is not on the cut locus of $J_0$
(Corollary~\ref{uncut}). 

\begin{thm}\label{ferm-bog}
Let $J_0,J_1\in\JJ_g$ and let $\gam=\{J_t\}_{0\le t\le1}$ be a geodesic
from $J_0$ to $J_1$, $t$ being proportional to the arc-length parameter.
Assume that the geodesic lies completely in the complement of the cut locus
of $J_0$.
Then\\
1. the parallel transport in $\HH$ along $\gam$ is 
$\UU_{J_1J_0}^\HH=\left(\det\hfjj\right)^{-1/4}\,\PP_{J_1J_0}$;\\
2. the parallel transport in $\sqrt{\KK^{-1}}$ along $\gam$ is
$\UU_{J_1J_0}^{\sqrt{\KK^{-1}}}=\sqrt{\det J_{1/2}/\ii\big|_{V_{J_0}^{1,0}}}$,
and 
\[ \bra\,\UU_{J_1J_0}^{\sqrt{\KK^{-1}}}\sqrt{\mu'_0},\sqrt{\mu_0}\,\ket
=\left(\det\hfjj\right)^{1/4}\,\bra\sqrt{\mu'_0},\sqrt{\mu_0}\,\ket  \]
for any $\mu_0,\mu'_0\in\medwedge^nV_{J_0}^{1,0}$;\\
3. the parallel transport in $\hat\HH$ along $\gam$, which is
$\UU_{J_1J_0}^{\hat\HH}=\UU_{J_1J_0}^\HH\otimes\UU_{J_1J_0}^{\sqrt{\KK^{-1}}}$,
corresponds to the pairing between $\hat\HH_{J_0}$ and $\hat\HH_{J_1}$ given by
\[ \bra\psi_1\otimes\sqrt{\mu_1},\psi_0\otimes\sqrt{\mu_0}\ket
   =\bra\psi_1,\psi_0\ket\bra\sqrt{\mu_1},\sqrt{\mu_0}\ket, \qquad \psi_l\in
   \HH_{J_l},\;\mu_l\in{\textstyle \medwedge^n}\,V_{J_l}^{1,0}\;\,(l=0,1). \]
\end{thm}

\noindent{\em Proof:}
1. Choosing a unitary basis $\{e_i\}$ of $V_{J_0}^{1,0}$, the geodesics in
$\JJ_g\cong\SO(2n)/\U(n)$ are given by Proposition~\ref{geod}.2.
We can assume without loss of generality (cf.~the proof of Theorem~3.4 in
\cite{KW}) that $n=2$ and $k=1$, $b=b_1>0$; the case $n=1$ is trivial.
Then as in the proof of Proposition~\ref{cut}, the vectors 
$e^{(t)}_1=\cos bt\;e_1-\sin bt\;\bar e_{\bar 2}$,
$e^{(t)}_2=\cos bt\;e_2+\sin bt\;\bar e_{\bar 1}$ form a unitary basis of
$V_{J_t}^{1,0}$.
Since $(\del P)\,e^{(t)}_1=-b\,e^{(t)}_2$ and
$(\del P)\,e^{(t)}_2=b\,e^{(t)}_1$, we have $(\del P)_1^{\;\bar 2}=-b$
and hence $(\del P)^{12}=b$; here the tensor indices correspond to the
basis $\{e^{(t)}_1,e^{(t)}_2\}$.
We want to find $\al(t)$ such that the quantity $\al(t)\bra\psi',\psi_t\ket$
is independent of $t$ for any $\psi'\in\HH_{J_0}$ if $\psi_t\in\HH_{J_t}$ is
a parallel transport of $\psi_0\in\HH_{J_0}$.
This would imply $\UU_{J_tJ_0}^\HH=\al(t)^{-1}\,\PP_{J_tJ_0}$.
Since $\psi_t$ satisfies the differential equation
\[ \dot\psi_t=-\hf(\del P)^{ij}\,\nabla_i\nabla_j\,\psi_t
             =-b\,\nabla_1\nabla_2\,\psi_t,     \]
we have
\begin{eqnarray*}
\frac{d}{dt}\big(\al(t)\bra\psi',\psi_t\ket\big)
\eq\dot\al(t)\bra\psi',\psi_t\ket-b\,\al(t)
                 \bra\psi',\nabla_1\nabla_2\,\psi_t\ket    \\
\eq\dot\al(t)\bra\psi',\psi_t\ket-b\,\al(t)\bra\psi',
   (\sec bt\,\nabla^{(0)}_1-\tan bt\,\nabla_{\bar2})\nabla_2\,\psi_t\ket   \\
\eq\dot\al(t)\bra\psi',\psi_t\ket+b\,\al(t)\tan bt\,
                 \bra\psi',\{\nabla_{\bar2},\nabla_2\}\psi_t\ket \\
\eq(\dot\al(t)-b\,\al(t)\tan bt)\,\bra\psi',\psi_t\ket.
\end{eqnarray*}
Solving $\dot\al-b\al\tan bt=0$ with $\al(0)=1$, we get
$\al(t)=(\cos bt)^{-1}$.
By Proposition~\ref{cut}, the assumption on the geodesic means that
$\frac{J_0+J_t}{2}$ is invertible for all $t\in[0,1]$.
Since $\det\frac{J_0+J_t}{2}=\cos^4bt$, this means $|b|<\frac{\pi}{2}$ and
the result follows.\\
2. The formula for $\UU_{J_1J_0}^{\sqrt{\KK^{-1}}}$ follows from
Lemma~\ref{A}.1.
It suffices to show the rest when $n=2$ using the above parametrisation.
If we take $\mu_0=e^{(0)}_1\wedge e^{(0)}_2$, then $\bra\mu_0,\mu_0\ket=1$ and
\[ \UU_{J_1J_0}^{\KK^{-1}}\,\mu_0=e^{(t)}_1\wedge e^{(t)}_2
        =\cos^2bt\,\mu_0+\cdots,\]
where the omitted part has a factor from $V_{J_0}^{0,1}$.
The result then follows from
\[  \bra\,\UU_{J_1J_0}^{\KK^{-1}}\,\mu_0,\mu_0\,\ket=\cos^2bt
               =\left(\det\textfrac{J_0+J_t}{2}\right)^{1/2}. \]
\noindent 3. This is an immediately consequence of parts~1 and 2.
\qed\vspace{10pt}

Notice that although power of the factor $\left(\det\hfjj\right)^{-1/4}$ in 
Theorem~\ref{ferm-bog}.1 is opposite to that in Theorem~\ref{bos-bog}.1,
both are greater than $1$.
Using the fermionic analog of the Bergman kernel projection in
Proposition~\ref{ferm-berg}, we have

\begin{cor}\label{ferm-UP}
Let $\psi=e^{\frac{\ii}{2}\varpi_{J_0}}\wedge\phi\in\HH_{J_0}$.
Then for fermionic but real $\tht\in\Pi V$,
\[ (\UU_{J_1J_0}^\HH\psi)(\tht)=\left(\det\hfjj\right)^{-1/4}\,
   e^{\frac{\ii}{4}\varpi_{J_1}(\tht)}\int_{\Pi V}\exp\big[\,
   g(\tht_{J_1}^{1,0},\chi)+\textfrac{\ii}{4}\varpi_{J_1}(\chi)
   +\textfrac{\ii}{4}\varpi_{J_0}(\chi)\,\big]\,\phi(\chi)
   \;\tilde\eps_g(\chi).                        \]
\end{cor}

We recall from \S\ref{B}.2 the notion of fermionic coherent states.

\begin{thm}\label{ferm-coh}
Under the assumptions of Theorem~\ref{ferm-bog},\\
1. the parallel transport along $\gam$ of the coherent state $c_{J_0}^\al$
is, for $\tht\in\Pi V$,
\[  (\UU^\HH_{J_1J_0}c_{J_0}^\al)(\tht)
    =\left(\det\hfjj\right)^{\!1/4}\,e^{\frac{\ii}{4}\varpi_{J_1}(\tht)}\,
    \exp\big[\textfrac{\ii}{2}g\big(\tht_{J_1}^{1,0}-\bar\al,\,
    \left(\hfjj\right)^{\!\!-1}\!(\tht_{J_1}^{1,0}-\bar\al)\,\big)\big];   \]
2. the parallel transport along $\gam$ of any state
$\psi=e^{\frac{1}{4}\varpi_{J_0}}\wedge\phi\in\HH_{J_0}$ is,
for $\tht\in\Pi V$,
\[ (\UU_{J_1J_0}^\HH\psi)(\tht)=\left(\det\hfjj\right)^{\!1/4}\!
    e^{\frac{1}{4}\varpi_{J_1}(\tht)}\!\!\int_{\Pi V}\!\!\!
    \exp\big[\textfrac{\ii}{2}g\big(\tht_{J_1}^{1,0}-\chi_{J_0}^{0,1},
    \left(\hfjj\right)^{\!\!-1}\!\!(\tht_{J_1}^{1,0}-\chi_{J_0}^{0,1})\big)
    +\textfrac{\ii}{4}\varpi_{J_0}(\chi)\big]\phi(\chi)
    \,\tilde\eps_g(\chi).                                         \]
\end{thm}

\noindent{\em Proof:}
1. We follow the proof of Theorem~\ref{bos-coh}.1.
By Lemma~\ref{pfA} and Proposition~\ref{ferm-berg}, we get, for fermionic but
real $\tht\in\Pi V$,
\begin{eqnarray*}
(\PP_{J_1J_0}c^\al_{J_0})(\tht)\eq e^{\frac{\ii}{4}\varpi_{J_1}(\tht)}
 \int_{\Pi V}\exp\big[\,g(\tht_{J_1}^{1,0},\chi)+\textfrac{\ii}{2}
 g\big(\hfjj\chi,\chi\big)\big]\,e^{g(\chi,\bar\al)}\;\tilde\eps_g(\chi) \\
\eq e^{\frac{\ii}{4}\varpi_{J_1}(\tht)}
 \exp\big[\textfrac{\ii}{2}g\big(\tht_{J_1}^{1,0}-\bar\al,\,
 \big(\hfjj\big)^{\!\!-1}\!(\tht_{J_1}^{1,0}-\bar\al)\big)\big]
 \int_{\Pi V}e^{\frac{\ii}{2}g\big(\hfjj\chi',\chi'\big)}\;\tilde\eps(\chi') \\
\eq\left(\det\hfjj\right)^{1/2}\,e^{\frac{\ii}{4}\varpi_{J_1}(\tht)}\,
 \exp\big[\textfrac{\ii}{2}g\big(\tht_{J_1}^{1,0}-\bar\al,\,
 \left(\hfjj\right)^{\!\!-1}\!(\tht_{J_1}^{1,0}-\bar\al)\,\big)\big],
\end{eqnarray*}
where we made a change of variable
$\chi'=\chi-\ii(\hfjj)^{-1}(\tht^{1,0}_{J_1}-\bar\al)$
and used Lemma~\ref{pfA}.
The condition that $\frac{J_0+J_t}{2}$ is invertible for all $0\le t\le1$
implies that $\hfjj$ is in the same connected component of invertible
skew-symmetric operators as $J_0$.
The result then follows from Theorem~\ref{ferm-bog}.1.\\
2. By Proposition~\ref{ferm-berg}, we have
\[  \psi(\tht,\bar\tht)=\int_{\Pi V}c_J^\chi(\tht)\,e^{-g(\chi,\bar\chi)}
    \,\phi(\chi,\bar\chi)\,\tilde\eps(\chi). \]
The result then follows from part~1 and the linearity of fermionic integration.
\qed\vspace{10pt}

When $n=2$, if the geodesic $\{J_t\}_{0\le t\le 1}$ from $J_0$ to $J_1$
is given by $z(t)=\tan bt$, where $|b|<\frac{\pi}{2}$, then
Theorem~\ref{ferm-coh} gives
\[  (\UU_{J_1J_0}c_{J_0}^\al)(\tht)=\cos b\; 
    \exp\big[(\tht^1\bar\al^{\bar1}+\tht^2\bar\al^{\bar2})\sec b
   +\hf(\bar\al^{\bar1}\bar\al^{\bar2}+\tht^1\tht^2)\tan b
   -\hf(\tht^1\bar\tht^{\bar1}+\tht^2\bar\tht^{\bar2})\big], \] 
which can also be obtained by solving the equation of parallel transport as
in the bosonic case (cf.~\cite{KW}).

\subsection{Fermionic systems with symmetries}

Let $(V,g)$ be a Euclidean space of dimension $2n$.
The action of $\SO(V,g)$ on $\JJ_g$ can be lifted to $\HH$, preserving the
connection.
It can be lifted to an action of the spin group $\Spin(V,g)$, which is a
double cover of $\SO(V,g)$, on $\sqrt{\KK^{-1}}$ and hence on $\hat\HH$.
The lifted action preserves the flat connection.

Suppose $K$ is a compact Lie group with Lie algebra $\kk$ acting on $(V,g)$
by an orthogonal representation.
Then it also acts on the bundles $\HH$ and $\hat\HH$ preserving the
connections.
Over the fixed-point set $(\JJ_g)^K$, the group $K$ acts on the fibres of
$\HH$.
Each fibre splits orthogonally into a direct sum of subspaces of various
representation-types of $K$.
Since $K$ preserves the connection, the restriction of the bundle $\HH$ to
$(\JJ_g)^K$ together with the projectively flat connection splits into
sub-bundles on which $K$ acts fibre-wisely.
In particular, we have a projectively flat sub-bundle $\HH^K\to(\JJ_g)^K$
on which $K$ acts trivially on the fibres.

We now study the fermionic reduced phase space $\Pi V\sq K$ and its
quantisation.
The action of the Lie algebra $\kk$ yields Hamiltonian ``vector fields"
on $\Pi V$ \cite{Ko77}.
The moment map $\mu_K$ is given by, for any $A\in\kk$,
$\bra\mu_K,A\ket=\frac{1}{2}g(A\cdot,\cdot)\in\medwedge^2V^*$ or
\[  \bra\mu_K(\tht),A\ket=\hf g(A\tht,\tht),\quad\tht\in\Pi V,\;A\in\kk  \]
using fermionic variables.
Following the construction of the usual symplectic quotients, the fermionic
analogue $\Pi V\sq K$ should be the ``spec" of the non-commutative ring
$(\medwedge^\bullet(V^*)^\CO/\bra\mu_K\ket)^K$, where $\bra\mu_K\ket$ is
the ideal generated by $\bra\mu_K,A\ket$ for all $A\in\kk$.
The ``space" $\Pi V\sq K$ is not usually a graded manifold in the sense of
\cite{Ko77}; it would be so if $0$ were a regular value of $\mu_K$ \cite{BS}.
So fermionic symplectic quotients are interesting examples of supermanifolds
with curved fermionic coordinates. 
Consider the example $V=\RE^{2n}$ with $K=S^1$ acting by weights
$\lam_1,\dots,\lam_r\ne0$, $\lam_{r+1}=\cdots=\lam_n=0$.
Then the ``coordinate ring" of $\Pi\RE^{2n}\sq S^1$ is generated by
$1,\tht^1\tht^2,\dots,\tht^{2r-1}\tht^{2r},\tht^{2r+1},\dots,\tht^{2n}$
subject to a relation $\lam_1\tht^1\tht^2+\cdots+\lam_r\tht^{2r-1}\tht^{2r}=0$.
Here $\tht^1,\dots,\tht^{2n}$ are the fermionic coordinates on $\Pi\RE^{2n}$.
When $r=1$, the above ring is the exterior algebra on $\tht^3,\dots,\tht^{2n}$
and thus $\Pi\RE^{2n}\sq S^1\cong\Pi\RE^{2n-2}$.

\begin{prop}
If $\dim(\JJ_g)^K>0$, then there is a non-zero $K$-invariant complex subspace
$(V',J_0)$ of $(V,J_0)$ with a $K$-invariant quaternionic structure such that
the restriction of $\mu_K$ to $\Pi V'$ is invariant under the scalar
multiplication by quaternions of unit norm. 
\end{prop}

\noindent{\em Proof:}
By Proposition~\ref{K-b-f}.2, there is a non-zero $K$-invariant complex
subspace $(V',J_0)$ of $(V,J_0)$ and a $K$-invariant quaternionic structure
$Q$ on $V'$ such that $g(Q\cdot,Q\cdot)=g(\cdot,\cdot)$ on $V'$.
For any $A\in\kk$, we have
\[  \bra\mu_K(Q\tht),A\ket=\hf g(Q\tht,AQ\tht)=\hf g(Q\tht,QA\tht)
    =\hf g(\tht,A\tht)=\bra\mu_K(\tht),A\ket,                         \]
where $\tht\in\Pi V'$. 
The result then follows easily from $g(J_0\cdot,J_0\cdot)=g(\cdot,\cdot)$
and $QJ_0=-J_0Q$.
\qed\vspace{5pt}

\section{Concluding remarks}\label{con}

We end with a comparison of the quantisation of bosons and fermions.
As explained in \S\ref{intro}, the existence of projectively flat connection
is due largely to the fact that the irreducible representation of the operator
algebra (Heisenberg algebra for bosons and Clifford algebra for fermions)
is unique up to unitary equivalence.
This enables us to identify, up to a phase, states in Hilbert spaces
constructed from various linear polarisations. 
Moreover, the geometric structure of the bundle of Hilbert spaces and results
on parallel transport are rather similar in the bosonic and fermionic cases.

We note however some conceptual differences.
For bosons, the positivity condition is on the polarisation whereas for
fermions, it is on the Euclidean structure $g$.
Indeed, the unitarity of the representation of the Heisenberg algebra is
not sensitive to the sign of $\om$, whereas for fermions, the positivity
condition on $g$ is the requirement for unitarity \cite{W91}.
On the other hand, the positivity condition on polarisation for bosons
guarantees the existence of holomorphic sections rather than elements
in higher cohomology groups.
For fermions, the cohomology is always at the zeroth degree; this reflects
the Dirac sea picture in physics.

Mathematically, the spaces of allowed polarisations are Hermitian symmetric
spaces in both cases.
For bosons, the space non-compact, non-positively curved.
Though it is contractible, the difficulty occurs at the boundary at infinity,
where the limit of parallel transport should be carefully taken \cite{KW}.
For fermions, the space of polarisation is compact, non-negatively curved.
Though it has no boundary, interesting phenomena (non-uniqueness of geodesics,
degeneracy of the half-form pairing) because of cut locus (\S\ref{A}.3 and
\S\ref{ferm}.3).
Furthermore, the half-form bundles in metaplectic correction are of opposite
powers of the canonical bundle in the two cases in order to cancel the
curvature of the projectively flat connection.
We hope that these observations are useful to future research on the
quantisation of more general symplectic or graded symplectic manifolds.

\vspace{30pt}
\setcounter{section}{0}
\renewcommand{\thesection}{\Alph{section}}

\noindent {\bf\Large Appendix}

\section{Geometry of the space of complex structures}\label{A}

\subsection{Complex structures on a vector space}

Let $V$ be a real vector space of dimension $2n$.
Consider the set $\JJ$ of complex structures on $V$ compatible with a given
orientation on $V$.
For each $J\in\JJ$, there is a decomposition $V^\CO=V_J^{1,0}\oplus V_J^{0,1}$,
where $V_J^{1,0}$, $V_J^{0,1}$ are the holomorphic and anti-holomorphic
subspaces, on which $J=\pm\ii$, respectively.
Similarly, there is a decomposition
$(V^*)^\CO=(V^{1,0}_J)^*\oplus(V^{0,1}_J)^*$.
For $x\in V$, $\al\in V^*$, we write, accordingly,
\[ x=x_J^{1,0}+x_J^{0,1},\quad\al=\al_J^{1,0}+\al_J^{0,1}.  \]

The space $\JJ$ is a connected smooth manifold of real dimension $2n^2$.
At any $J_0\in\JJ$, the tangent space of $\JJ$ is
$T_{J_0}\JJ\cong \Hom_\CO(V_{J_0}^{1,0},V_{J_0}^{0,1})$.
Moreover, a dense open subset of $\JJ$ can be parametrised by
$Z\in\Hom_\CO(V_{J_0}^{1,0},V_{J_0}^{0,1})$: for any such $Z$, the
corresponding subspace $V_J^{1,0}$ is the graph of $Z$, that is, $V_J^{1,0}$
consists of the vectors of the form ${x\choose Zx}$ under the decomposition
$V^\CO=V_{J_0}^{1,0}\oplus V_{J_0}^{0,1}$, where $x\in V_{J_0}^{1,0}$.
These open sets (for various $J_0$) have complex coordinates and form an open
cover of $\JJ$.
This defines a complex structure on $\JJ$.
For a fixed $J_0\in\JJ$, the complement of the open set consists of $J\in\JJ$
such that $V_J^{1,0}\cap V_{J_0}^{0,1}\ne\{0\}$, which happens when $J+J_0$
is not invertible.
On the other hand, not every $Z\in\Hom_\CO(V_{J_0}^{1,0},V_{J_0}^{0,1})$
corresponds to a complex structure. 
If it does, then the condition $V_J^{1,0}\cap V_J^{0,1}=\{0\}$ implies that
$1-\bar ZZ\in\End(V_{J_0}^{1,0})$ is invertible.
An element $Z$ that violates this condition is on the ``boundary'' of $\JJ$.
Finally, the projection onto $V_J^{1,0}$ along $V_J^{0,1}$ is
$P_J=\hf(1-\ii J)$.
Given $J_0$, on the dense set where $J$ can be parametrised by $Z$,
the projection is
\[ P_J={1\choose Z}(1-\bar ZZ)^{-1}(1,\;\,-\bar Z).\]

Suppose $\del J$ is an infinitesimal variation of $J\in\JJ$.
(The symbol $\del$ can be interpreted as the differential on $\JJ$.)
Since the change $\del P=-\frac{\ii}{2}\del J$ of $P=P_J$ anti-commutes
with $J$, it is off-diagonal with respect to the decomposition
$V^\CO=V_J^{1,0}\oplus V_J^{0,1}$.
The new holomorphic subspace $V_{J+\del J}^{1,0}$ is the graph of the
$\Hom_\CO(V_J^{1,0},V_J^{0,1})$-component of $\del P$.
Since $P+\bar P=\id_V$, we have $\del\bar P=-\del P$.

We consider the vector bundle $\VV\to\JJ$ whose fibre over $J\in\JJ$ is
$V_J^{1,0}$.
This is a sub-bundle of the product bundle $\JJ\times V^\CO$ and has a
connection defined by the projection $P_J$.
This connection on $\VV$ is $\A^\VV=-(\del P_J)P_J$ and its curvature is
\begin{eqnarray*}
\F^\VV\eq P_J\,\del P_J\wedge\del P_J\,P_J=-\qt P_J\,\del J\wedge\del J\,P_J \\
         \eq-{1\choose Z}(1-\bar ZZ)^{-1}\del\bar Z
             \wedge(1-Z\bar Z)^{-1}\del Z(1-\bar ZZ)^{-1}(1,\;\,-\bar Z).
\end{eqnarray*}
The curvature of the line bundle $\det\VV=\medwedge^n\VV\to\JJ$ is the $2$-form
\[ \F^{\det\!\VV}=\tr(P_J\,\del P_J\wedge\del P_J\,P_J)=
    \tr((1-Z\bar Z)^{-1}\del Z\wedge(1-\bar ZZ)^{-1}\del\bar Z).        \]
(All expressions in terms of $Z$ are valid on a dense open set of $\JJ$ only.)
The canonical line bundle over $V$ is the product bundle $V\times\KK_J$, where
$\KK_J=\medwedge^n(V^{1,0}_J)^*$.
We have a line bundle $\KK\to\JJ$ whose fibre over $J\in\JJ$ is $\KK_J$.
In fact, $\KK=(\det\!\VV)^*$ and its curvature is $\F^\KK=-\F^{\det\!\VV}$.

The space $\JJ$ has a transitive action of $\GL_+(V)$, the identity component
of $\GL(V)\cong\GL(2n,\RE)$ preserving the orientation on $V$.
At $J_0\in\JJ$, the isotropy subgroup $\GL(V,J_0)\cong\GL(n,\CO)$ consists of
elements in $\GL(V)$ that commutes with $J_0$.  
Therefore $\JJ$ can be identified as the homogeneous space
$\GL_+(V)/\GL(V,J_0)$.
The Lie algebra of $\GL_+(V)$ has the decomposition
$\gl(V)=\gl(V,J_0)\oplus\mm$, where $\gl(V,J_0)$, $\mm$ are the subspaces in
$\gl(V)$ of elements that commute, anti-commute with $J_0$, respectively.
Since $\mm$ is invariant under the adjoint action of $\GL(V,J_0)$ and since
$[\mm,\mm]\subset\gl(V,J_0)$, $\GL_+(V)/\GL(V,J_0)$ is a reductive symmetric
space.
The trace form on $\mm\cong T_{J_0}\JJ$ is a $\GL(V,J_0)$-invariant
non-degenerate symmetric bilinear form and it induces a pseudo-Riemannian
metric
\[   \upeta=2\tr(P_J\,\del P_J\,\del P_J\,P_J)
               =2\tr((1-Z\bar Z)^{-1}\del Z\,(1-\bar ZZ)^{-1}\del\bar Z)   \]
on $\JJ$.
It is pseudo-K\"ahler with the pseudo-K\"ahler form
\[ \upsigma=\F^\KK/\ii=\ii\tr(P_J\,\del P_J\wedge\del P_J\,P_J)
        =\ii\tr((1-Z\bar Z)^{-1}\del Z\wedge(1-\bar ZZ)^{-1}\del\bar Z).   \]
The group $\GL_+(V)$ acts on $Z$ by fractional linear transformations
preserving $\upsigma$.

The map $J\mapsto J_0^{-1}JJ_0$ is an isometric reflection on $\JJ$ that
induces minus the identity map on the tangent space $T_{J_0}\JJ\cong\mm$.
As $\JJ$ is reductive, a geodesic in $\JJ$ is of the form $t\mapsto[g_t]$ with 
$g_t=g_0\,e^{tM}$ (see \cite{No}), where $g_0\in\GL_+(V)$, $M\in\mm$ and $[g]$
denotes the coset in $\GL_+(V)/\GL(V,J_0)\cong\JJ$ represented by
$g\in\GL_+(V)$.
The parameter $t\in\RE$ is proportional to the arc-length on the geodesic.

\begin{prop}\label{J/2}
For any $t$, $J_{t/2}/\ii$ is the parallel transport in $\VV$ along the
geodesic from $J_0$ to $J_t$.
\end{prop}

\noindent{\em Proof:}
Since the geodesic starts from $J_0$, we have $g_0=1$, hence
$g_tg_{t'}=g_{t+t'}$ and $g_t^{-1}=g_{-t}$.
The reflection at $J_0$ reverses the direction of each geodesic passing
through $J_0$, as $J_0^{-1}g_tJ_0=g_{-t}$.
The complex structure at $[g_t]$ is $J_t=g_tJ_0g_t^{-1}=g_{2t}J_0=J_0g_{-2t}$,
and hence $J_0J_{-t}=J_tJ_0=-g_{2t}$.
$J_0/\ii$ is the identity map on $V_{J_0}^{1,0}$, $J_{t/2}=g_tJ_0$ maps
$V_{J_0}^{1,0}$ to $V_{J_t}^{1,0}$, and $\frac{d}{dt}J_{t/2}=g_tMJ_0$ maps
$V_{J_0}^{1,0}$ to $V_{J_t}^{0,1}$.
So $J_{t/2}/\ii$ satisfies all the requirements that uniquely defines the
parallel transport in $\VV$.
\qed\vspace{5pt}

\begin{cor}
If $J$, $J'$ and $J''$ are three points on a geodesic such that $J'$
bisects the segment between $J$ and $J''$, then $J'J=J''J'$ and $J'/\ii$,
which maps $V_J^{1,0}$ to $V_{J''}^{1,0}$, is the parallel transport
in $\VV$ from $J$ to $J''$ along the geodesic.
\end{cor}

\subsection{Complex structures compatible with a symplectic 
or Euclidean structure}

Given a symplectic form $\om$ on $V$, a complex structure $J$ on $V$ is
compatible to $\om$ if $\om(J\cdot,J\cdot)=\om(\cdot,\cdot)$ and
$\om(\cdot,J\cdot)>0$.
(The second condition implies that the orientation defined by $J$ coincides 
with that of the volume form $\om^n/n!$.)
Let $\JJ_\om$ be the set of such $J$.
The symplectic form defines an isomorphism $\nu=\nu_\om\colon V\to V^*$ by
$x\in V\mapsto\iota_x\om\in V^*$ and a holomorphic involution $\s=\s_\om$
on $\JJ$ by $J\mapsto\nu^{-1}\circ{}^T\!J\circ\nu$. 
The fixed-point set $\JJ^\s$ consists of $J\in\JJ$ satisfying
$\om(J\cdot,J\cdot)=\om(\cdot,\cdot)$.
The connected components of $\JJ^\s$, one of which is $\JJ_\om$, are labelled
by the signature of the symmetric bilinear form $\om(\cdot,J\cdot)$.
Since $\s$ is an isometry, $\JJ^\s$ is a totally geodesic submanifold.
In fact, $\JJ_\om$ is a K\"ahler manifold since the pseudo-K\"ahler metric
$\upeta$ restricts to a positive-definite metric $\upeta_\om$ on $\JJ_\om$.
So the restriction $\upsigma_\om$ of $\upsigma$ to $\JJ_\om$ is a K\"ahler
form.
Given $J_0\in\JJ_\om$, the whole space $\JJ_\om$ can be parametrised by
$Z\in\Hom_\CO(V_{J_0}^{1,0},V_{J_0}^{0,1})$.
We denote by the same notation $\nu\colon V_{J_0}^{0,1}\to(V_{J_0}^{1,0})^*$
the restriction of the isomorphism $\nu\colon V^\CO\to(V^*)^\CO$. 
Then $J\in\JJ^\s$ if and only if $\nu\circ Z\in\sym^2(V_{J_0}^{1,0})^*$, or
equivalently, $Z\circ\nu^{-1}\in\sym^2V_{J_0}^{0,1}$.
The tangent space $T_{J_0}\JJ_\om$ consists of $Z$ satisfying this condition.
The condition $\om(\cdot,J\cdot)>0$ is equivalent to $1-\bar ZZ>0$ with respect
to the Hermitian form $h_0(x,y)=\om(x,J_0\bar y)=-\ii\om(x,\bar y)$,
$x,y\in V_{J_0}^{1,0}$.
The symplectic group $\Sp(V,\om)$ acts transitively on $\JJ_\om$ and the
isotropic subgroup at $J_0$ is the unitary group $\U(V,h_0)$.
So $\JJ_\om\cong\Sp(V,\om)/\U(V,h_0)$.
It can be identified holomorphically with a bounded Hermitian symmetric domain
or with the Siegel upper-half space; the two are related by a Cayley transform.

If instead there is a Euclidean inner product $g$ on $V$, let $\JJ_g$ be
the set of complex structures $J$ that is compatible with $g$, i.e.,
$g(J\cdot,J\cdot)=g(\cdot,\cdot)$, and the orientation on $V$.
There is an isomorphism $\nu=\nu_g\colon V\to V^*$ defined
by $x\in V\mapsto\iota_xg\in V^*$ and an involution
$\s=\s_g\colon J\mapsto\nu^{-1}\circ{}^T\!J\circ\nu$ on $\JJ$.
The fixed-point set is $\JJ^\s=\JJ_g$.
As $\s$ is an isometry, $\JJ_g$ is totally geodesic in $\JJ$ as in the
symplectic case.
$\JJ_g$ is K\"ahler since the restriction $\upeta_g$ of $-\upeta$ to $\JJ_g$
is positive definite; the restriction $\upsigma_g$ of $-\upsigma$ to $\JJ_g$
is the K\"ahler form.
Given $J_0\in\JJ_g$, we denote also by
$\nu\colon V_{J_0}^{0,1}\to(V_{J_0}^{1,0})^*$
the restriction of $\nu\colon V^\CO\to(V^*)^\CO$.
On the dense set of $\JJ$ that can be parametrised by $Z$, $J\in\JJ_g$
if and only if the corresponding $Z$ satisfies
$\nu\circ Z\in\medwedge^2(V_{J_0}^{1,0})^*$, or equivalently,
$Z\circ\nu^{-1}\in\medwedge^2V_{J_0}^{0,1}$.
The tangent space $T_{J_0}\JJ_g$ consists of $Z$ satisfying this condition.
For any such $Z$, we always have $1-\bar ZZ>0$ with respect to the Hermitian
form $h_0(x,y)=g(x,\bar y)$, $x,y\in V_{J_0}^{1,0}$.
The group $\SO(V,g)$ acts transitively on $\JJ_g$ and the isotropic subgroup
at $J_0$ is $\U(V,h_0)$.
The space $\JJ_g\cong\SO(V,g)/\U(V,h_0)$ is a compact Hermitian symmetric
space.
Finally, if in addition there is a symplectic form $\om$ on $V$ such that
$\om(\cdot,J_0\cdot)$ is proportional to $g$, then $\JJ_\om$ and $\JJ_g$
intersects at $J_0$ orthogonally with respect to the pseudo-K\"ahler
metric $\upeta$ on $\JJ$.

We describe the results using tensor indices.
Let $\{e_i\}_{1\le i\le n}$ be a basis of $V_{J_0}^{1,0}$.
Then $\{\bar e_{\bar i}\}$ is a basis of $V_{J_0}^{0,1}$.
We represent $Z\in\Hom_\CO(V_{J_0}^{1,0},V_{J_0}^{0,1})$ by a matrix
$Z_i^{\;\bar j}$ such that $Ze_i=Z_i^{\;\bar j}\bar e_{\bar j}$.
If $(V,\om)$ is a symplectic vector space, we have
$\om_{\bar ij}=\om(\bar e_{\bar i},e_j)=-\om_{j\bar i}$.
Set $Z_{ij}=(\nu\circ Z)_{ij}=Z_i^{\;\bar k}\om_{\bar kj}$ and
$Z^{\bar i\bar j}=(Z\circ\nu^{-1})^{\bar i\bar j}
=\om^{\bar ik}Z_k^{\;\,\bar j}$.
Then $Z$ determines an element in $\JJ_\om$ if and only if $Z_{ij}=Z_{ji}$
(or $Z^{\bar i\bar j}=Z^{\bar j\bar i}$) and the matrix 
$\del_i^{\;j}-Z_i^{\;\bar k}\bar Z_{\bar k}^{\;j}$ is positive definite.
If $(V,g)$ is a Euclidean space instead, then
$g_{\bar ij}=g(\bar e_{\bar i},e_j)=g_{j\bar i}$.
Set $Z_{ij}=(\nu\circ Z)_{ij}=Z_i^{\;\bar k}g_{\bar kj}$ and
$Z^{\bar i\bar j}=(Z\circ\nu^{-1})^{\bar i\bar j}=g^{\bar ik}Z_k^{\;\,\bar j}$.
Then $Z$ determines an element in $\JJ_g$ if and only if $Z_{ij}=-Z_{ji}$
(or $Z^{\bar i\bar j}=-Z^{\bar j\bar i}$).
If there is a variation $\del J$ of $J\in\JJ$, then we have tensors
$(\del P)_i^{\;\bar j}=-(\del\bar P)_i^{\;\bar j}$ and
$(\del P)_{\bar i}^{\;j}=-(\del\bar P)_{\bar i}^{\;j}$.
We note that $\{e_i+\del e_i\}$, where
$\del e_i=(\del P)_i^{\;\bar j}\bar e_{\bar j}$, is a basis of
the new holomorphic subspace $V_{J+\del J}^{1,0}$, whereas
$\{\bar e_{\bar i}+\del\bar e_{\bar i}\}$, where
$\del\bar e_{\bar i}=(\del\bar P)_{\bar i}^{\;\,\bar j}\bar e_{\bar j}
=-(\del P)_{\bar i}^{\;\,\bar j}\bar e_{\bar j}$,
is a basis of $V_{J+\del J}^{0,1}$. 
If $(V,\om)$ is symplectic and $J\in\JJ_\om$, then $J+\del J\in\JJ_\om$ 
(to the first order) if and only if any of the tensors
$(\del P)_{ij},(\del P)^{ij},(\del P)^{\bar i\bar j},(\del P)_{\bar i\bar j}$
is symmetric.
If $(V,g)$ is Euclidean and $J\in\JJ_g$, then $J+\del J\in\JJ_g$ (to the
first order) if and only if any of the above tensors is anti-symmetric.

Choosing a unitary basis $\{e_i\}_{1\le i\le n}$ of $V_{J_0}^{1,0}$ in both the
symplectic and the orthogonal cases, we have $\JJ_\om\cong\Sp(2n,\RE)/\U(n)$
and $\JJ_g\cong\SO(2n)/\U(n)$, respectively, where $J_0$ is identified with the
coset $o$ of the identity element.
Using the basis $\{e_i,\bar e_{\bar i}\}$ of $V^\CO$, the Lie groups and/or
their Lie algebras that appear in the above identifications are
\[ \U(n)=\bigset{\mat{U}{0}{0}{\bar U}}{\,{}^T\!\bar UU=I_n},
\quad \uu(n)=\bigset{\mat{A}{0}{0}{\bar A}}{\,{}^T\!\bar A=-A}, \]
\[ \rsp(2n,\RE)=\bigset{\mat{A}{B}{\bar B}{\bar A}}{\,
   \substack{{}^T\!\bar A=-A,\vspace{3pt} \\{}^T\!B=B}},  \quad
\so(2n)=\bigset{\mat{A}{B}{\bar B}{\bar A}}{\,\substack{{}^T\!\bar A=-A,
                      \vspace{3pt}\\{}^T\!B=-B}}, \]
where $U,A,B$ are $n\times n$ complex matrices.
The following results on geodesics are well-known:

\begin{prop}\label{geod} {\rm (\cite{Si,H})}
1. A geodesic $\gam$ in $\Sp(2n,\RE)/\U(n)$ from $o$ is of the form
\[ \gam(t)=\Big[k\,\mat{\cosh Bt}{\sinh Bt}{\sinh Bt}{\cosh Bt}k^{-1}\Big], \]
where $k\in\U(n)$ and $B=\diag\{b_1,\cdots,b_r,0,\cdots,0\}$ for some
$b_1,\dots,b_r>0$, $r\le n$.\\
2. A geodesic $\gam$ in $\SO(2n)/\U(n)$ from $o$ is of the form
\[  \gam(t)=\Big[k\,\Big(\substack{\quad\cos\sqrt{-B^2}t\qquad
    \frac{B}{\sqrt{-B^2}}\sin\sqrt{-B^2}t\!\!\!\\ \!\!\!\!\!
    \frac{B}{\sqrt{-B^2}}\sin\sqrt{-B^2}t\qquad\cos\sqrt{-B^2}t}
    \;\;\Big)k^{-1}\Big], \]
where $k\in\U(n)$ and $B=\diag\Big\{\mat{}{\;\,b_1\!\!\!}{-b_1}{},
\cdots,\mat{}{\;\,b_r\!\!\!}{-b_r}{},0,\cdots,0\,\Big\}$ for some
$b_1,\dots,b_r>0$, $r\le[n/2]$.
(In this case, $\sqrt{-B^2}=\diag\{b_1,b_1,\cdots,b_r,b_r,0,\cdots,0\}$.)
\end{prop}

\noindent{\em Proof:}
Writing $\rsp(2n,\RE)=\uu(n)\oplus\mm$ and $\so(2n)=\uu(n)\oplus\mm$,
respectively, geodesics are of the form $\gam(t)=[e^{tM}]$ for some
$M=\mat{0}{B'}{\!\overline{B'}}{\!0}\in\mm$.
Since ${}^T\!B'=\pm B'$, by Theorems~5 and 7 in \cite{H}, respectively,
there exists an $n\times n$ complex matrix $U$, ${}^T\bar UU=I_n$, such
that $B'=UB\,{}^T\!U$, where $B$ is of the required form.
The results then follow from simple calculations with
$k=\mat{U}{0}{0}{\bar U}\in\U(n)$.
\qed\vspace{5pt}

When $n=1$, $\JJ_\om=\JJ$ because every complex structure compatible with the
orientation is compatible with the symplectic form.
Choosing a base vector $e_1$ of the $1$-dimensional vector space
$V_{J_0}^{1,0}$, $\JJ_\om$ can be parametrised by $z=Z_1^{\;\bar1}\in\CO$
such that $|z|<1$, with $z=0$ for $J_0$.
The K\"ahler form and metric are, respectively,
\[ \upsigma_\om=\frac{2\ii dz\wedge d\bar z}{(1-|z|^2)^2},\qquad
   \upeta_\om=\frac{4\,dz\,d\bar z}{(1-|z|^2)^2}.   \]
A geodesic through $z=0$ is of the form $z(t)=e^{\ii\al}\tanh t$
($0\le\al<2\pi$), where $t\in\RE$ is half of the arc-length parameter.
For Euclidean space $(V,g)$, the first non-trivial case is, when $n=2$,
$\JJ_\om\cong\SO(4)/\U(2)=S^2$.
Choose a basis $\{e_1,e_2\}$ of $V_{J_0}^{1,0}$ such that
$g_{1\bar1}=g_{2\bar2}$, $g_{1\bar2}=g_{2\bar1}=0$.
Then the dense subset $\JJ_\om\backslash\{-J_0\}$ can be parametrised by
$Z=\mat{}{\;\,z\!\!\!}{-z}{}$, where $z\in\CO$.
On $\JJ_\om$, the point $-J_0$ (which would be $z=\infty$) is conjugate
to $J_0$ ($z=0$).
The K\"ahler form and metric are, respectively,
\[ \upsigma_g=\frac{2\ii dz\wedge d\bar z}{(1+|z|^2)^2},\qquad
   \upeta_\om=\frac{4\,dz\,d\bar z}{(1+|z|^2)^2}.   \]
A geodesic through $z=0$ is of the form $z(t)=e^{\ii\al}\tan t$
($0\le\al<2\pi$), where $t\in\RE$ is half of the arc-length parameter.
Note that $z(t)=\infty$ at $t=\frac{\pi}{2}$ corresponds to the antipodal
point $-J_0$ of $J_0$.

\subsection{Cut and first conjugate loci in $\JJ_g$}

Given a point $o$ in a Riemannian manifold $M$, the first conjugate point $p$
of $o$ along a geodesic $\gam$ from $o$ is a point such that there is a Jacobi
field along $\gam$ that is zero at $o$ and $p$ but nowhere zero in between.
The collection of such points form the first conjugate locus of $o$.
The cut point of $o$ along a geodesic $\gam$ from $o$ is the point $p$ such
that $\gam$ is length-minimising between $o$ and $p$ but fails to be so
beyond $p$.
There is an open cell $B$ in $T_oM$ such that the exponential map is a
diffeomorphism from $B$ onto a (connected) open subset of $M$ whose
compliment is the cut locus of $o$.
The image of the closure $\bar B$ under the exponential map is $M$.

While the structure of cut loci or first conjugate loci for general Riemannian
manifolds is quite complicated (see for example \cite{We}), there is a
Lie-theoretical description for compact Riemannian symmetric spaces.
For simply connected symmetric spaces (such as the space $\JJ_g$ above),
the cut locus and first conjugate locus coincide \cite{Cri}, though this fails
to be true in general \cite{Sa}.
For example, the cut and first conjugate loci of Grassmannian manifolds
are known explicitly in terms of Schubert varieties \cite{Wo} (see however
Remark~4.3 of \cite{Sa}).
We determine the cut (or the first conjugate) locus of $\JJ_g$, which
is the space of polarisations of fermionic systems.

\begin{prop}\label{cut}
Let $J_0,J\in\JJ_g$.
The following statements are equivalent:

(a) $J$ is on the cut locus of $J_0$;

(b) $\det\big(\frac{J_0+J}{2}\,\big)=0$;

(c) the pairing between $\KK_{J_0}^{-1}$ and $\KK_J^{-1}$ is degenerate.
\end{prop}

\noindent{\em Proof:}
Choosing a unitary basis $\{e_i\}_{1\le i\le n}$ of $V_{J_0}^{1,0}$, we
have $\JJ_g\cong\SO(2n)/\U(n)$.
The geodesics from $o$ are given by Proposition~\ref{geod}.2.
Since $k\in\U(n)$ acts as an isometry, we can assume $k=1$ as well as
$b_1\ge\cdots\ge b_r>0$ without loss of generality.
For any $t$, $V_{J_t}^{1,0}$ has a unitary basis consisting of vectors
$e^{(t)}_{2i-1}=\cos b_it\;e_{2i-1}-\sin b_it\;\bar e_{\overline{2i}}$,
$e^{(t)}_{2i}=\cos b_it\;e_{2i}+\sin b_it\;\bar e_{\overline{2i-1}}$
($1\le i\le r$) and $e^{(t)}_j=e_j$ ($2r+1\le j\le n$).

(a)$\Rightarrow$(b): 
The cut point of $o$ along the above geodesic is at $t=\pi/2b_1$.
It is clear that
$e^{(\pi/2b_1)}_1=-\bar e_{\bar2}\in V_{J_{\pi/2b_1}}^{1,0}\cap V_{J_0}^{0,1}$
and hence $\det\big(\frac{J_0+J_{\pi/2b_1}}{2}\,\big)=0$.

(b)$\Rightarrow$(a):
Consider a geodesic $\gam$ from $J_0$ to $J$ of the above form.
Then $\det\big(\frac{J_0+J}{2}\,\big)=0$ implies that $t=\pi/2b_i$ for some
$i=1,\dots,r$.
Assume $i=1$.
For $2\le j\le r$, let $b'_j$ be defined such that $|b'_j|\le b_1$ and
$b'_j=b_j\!\!\mod2b_1$.
Let $\gam'$ be the geodesic from $o$ corresponding to
$B'=\diag\Big\{\mat{}{\;\,b_1\!\!\!}{-b_1}{},\mat{}{\;\,b'_2\!\!\!}{-b'_2}{},
\cdots,\mat{}{\;\,b'_r\!\!\!}{-b'_r}{},0,\cdots,0\,\Big\}$.
Then $J=\gam'(\pi/2b_1)$ is the cut point of $o$ along $\gam'$.

(b)$\Leftrightarrow$(c):
Along the geodesic, let $\mu_t=e_1^{(t)}\wedge\cdots\wedge e_n^{(t)}$.
Then
\[ \bra\mu,\mu_0\ket=\prod_{i=1}^r\cos^2b_i
                    =\det\big(\textfrac{J_0+J}{2}\,\big)      \]
and hence the result.
\qed\vspace{5pt}

\begin{cor}\label{uncut}
If $J\in\JJ_g$ is not on the cut locus of $J_0$, then\\
1. $\det\big(\frac{J_0+J}{2}\,\big)>0$;\\
2. the inner product on $\sqrt{\KK_{J_0}^{-1}}$ extends continuously to a
non-degenerate pairing between $\sqrt{\KK_{J_0}^{-1}}$ and $\sqrt{\KK_J^{-1}}$.
\end{cor}

\noindent{\em Proof:}
Consider the geodesic in the proof of Proposition~\ref{cut}.\\
1. $\det\big(\frac{J_0+J}{2}\,\big)=\prod_{i=1}^r\cos^2b_i>0$ if $t<\pi/2b_i$
for all $i=1,\dots,r$.\\
2. Consider $\mu_t$ in the proof of Proposition~\ref{cut}.
The pairing is given by
$\bra\sqrt{\mu_t},\sqrt{\mu_0}\ket=\prod_{i=1}^r\cos b_i$.
\qed\vspace{5pt}

\section{Berezin integral and the fermionic Bergman kernel}\label{B}

\subsection{Calculus of fermionic variables}

Let $V$ be an $n$-dimensional real vector space with a (non-zero) volume
element $\eps=\eps_V\in\medwedge^n V$.
The Berezin integral of a form $\al\in\medwedge^\bullet(V^*)^\CO$ on $V$ is
\[   \int_{\Pi V}\al\,\eps_V=\bra\al^{(n)},\eps_V\ket,     \]
where the pairing is between the top-degree component $\al^{(n)}$ and $\eps_V$.
To highlight its formal similarity with the usual integration, the Berezin
integral is often expressed, as in the physics literature, as an
``integration" over fermionic variables.
While the setting is well known, we recall it here to fix the sign convention.
For a standard reference, see for example, \S1.4-7 of \cite{ZJ}.
For a mathematical treatment of graded manifolds or supermanifolds,
especially in the context of geometric quantisation, see \cite{Ko77}.

We imagine a copy $\Pi V$ of the vector space that is identical as $V$ except
it has fermionic coordinates, which are ``numbers" satisfying the same law of
addition but anti-commute when they are multiplied.
If we choose a basis $\{e_i\}_{1\le i\le n}$ of $V$, then a ``vector"
$\tht\in\Pi V$ has the form $\tht=\tht^ie_i$, where $\tht^1,\dots,\tht^n$
the fermionic coordinates.
Although $\Pi V$ does not exist as a set of points, the ``functions"
on $\Pi V$ are elements of the exterior algebra $\medwedge^\bullet(V^*)^\CO$.
In fact, any form
\[  \al=\sum_{k=0}^n\sum_{1\le i_1<\cdots<i_k\le n}\al_{i_1\dots i_k}
          e^{i_1}\wedge\cdots\wedge e^{i_k}   \]
on $V$ determines a ``function"
\[  \al(\tht)=\sum_{k=0}^n\sum_{1\le i_1<\cdots<i_k\le n}\al_{i_1\dots i_k}
          \tht^{i_1}\cdots\tht^{i_k}  \]
on $\Pi V$.
The ``derivative" of such functions corresponds the usual contraction on forms:
\[    \frac{\partial}{\partial\tht^i}\al(\tht)=(\iota_{e_i}\al)(\tht).   \]
Suppose the basis $\{e_i\}$ spans a unit volume, i.e.,
$\eps=e_1\wedge\cdots\wedge e_n$.
The volume element provides a ``measure" $\eps(\tht)=d\tht^1\cdots d\tht^n$ on
$\Pi V$.
The fermionic integral is defined by
\[  \int_{\Pi V}\al(\tht)\,\eps(\tht)
   =\int_{\Pi\RE^n}\al_{12\dots n}\tht^1\cdots\tht^n\,d\tht^1\cdots d\tht^n
   =(-1)^{\frac{n(n-1)}{2}}\al_{12\dots n}.    \]
This differs from the Berezin integral by a sign because $d\tht^i$ also
anti-commutes with $\tht^j$.  
As a useful example, we calculate
\[ \int_{\Pi\RE^2}e^{\ii a\tht^1\tht^2}\ii d\tht^1 d\tht^2
=\int_{\Pi\RE^2}(1+\ii a\tht^1\tht^2)\ii d\tht^1 d\tht^2
=a\int_{\Pi\RE}\tht^1\,d\tht^1\int_{\Pi\RE}\tht^2\,d\tht^2=a,  \]
where $a\in\RE$.

Suppose $V$ is even dimensional, say $\dim V=2n$.
Let $g$ be a Euclidean metric on $V$ and $\eps=\eps_g$, a unit volume element.
Set $\tilde\eps_g=\ii^n\,\eps_g$.

\begin{lem}\label{pfA}
If $A\in\End(V)$ is skew-symmetric with respect to $g$, then 
\[ \int_{\Pi V}e^{\frac{\ii}{2}g(A\tht,\tht)}\,\tilde\eps_g(\tht)=\Pf(A).  \]
\end{lem}

\noindent{\em Proof:}
We choose an orthonormal basis of $V$ so that $A$ decomposes as a direct sum
of $2\times2$ skew-symmetric matrices.
The result then follows from the example computed above.
\qed\vspace{10pt}

Now assume that $J$ is a complex structure on $V$ compatible with $g$ and
the orientation given by $\eps_g$.
If $A$ is invertible and if $A$ and $J$ are in the same connected component
of invertible, skew-symmetric operators on $V$, then
\[   \Pf(A)=(\det A)^{1/2},   \]
where the square root is chosen so that $(\det J)^{1/2}=1$.
Compare this with the usual Gaussian integral in the proof of
Theorem~\ref{bos-coh}.1, in which $A$ is symmetric and the determinant
factor $(\det A)^{1/2}$ appears in the denominator.

\subsection{The fermionic Bergman kernel and projection}

We work with the pre-quantum data of the fermionic system in \S\ref{ferm}.1:
a real Euclidean space $(V,g)$ of dimension $2n$ with a complex structure $J$
compatible with $g$ and a unit volume element $\eps_g$ which agrees with the
orientation of $J$.
Choosing a basis $\{e_i\}_{1\le i\le n}$ of $V_J^{1,0}$, we have complex
fermionic coordinates $\tht^1,\dots,\tht^n$ of $\tht\in\Pi V_J^{1,0}$ or
$\Pi V$.
An element $\psi$ in the pre-quantum Hilbert space $\HH_0$ can be regarded
as a ``function" $\psi(\tht,\bar\tht)$ of $\tht^i$ and $\bar\tht^{\bar i}$
($1\le i\le n$).
The covariant derivative along $e_i$ and $\bar e_{\bar j}$ are, respectively,
\[  \nabla_i=\textfrac{\pdr}{\pdr\tht^i}
             -\hf g_{i\bar j}\bar\tht^{\bar j},\qquad\nabla_{\bar j}
  =\textfrac{\pdr}{\pdr\bar\tht^{\bar j}}-\hf g_{i\bar j}\tht^i,   \]
where $g(\tht,\bar\tht)=g_{i\bar j}\tht^i\bar\tht^{\bar j}$.
Any $\psi=e^{\frac{\ii}{2}\varpi_J}\wedge\phi\in\HH_J$ can be written as
(cf.\ Theorem~\ref{ferm-holo}.1)
\[ \psi(\tht,\bar\tht)=\phi(\tht)\,e^{-\frac{1}{2}g(\tht,\bar\tht)},  \]
where $\phi(\tht)$ is a ``holomorphic function", that is, it depends on
$\tht^i$ only.
By Theorem~\ref{ferm-holo}.2, the inner product of
$\psi=e^{\frac{\ii}{2}\varpi_J}\wedge\phi$
and $\psi'=e^{\frac{\ii}{2}\varpi_J}\wedge\phi'$ in $\HH_J$ is
\[ \bra\psi,\psi'\ket=\int_{\Pi V}\phi(\tht)^*\,\phi'(\tht)
                      \;e^{-g(\tht,\bar\tht)}\,\tilde\eps_g(\tht).   \]
Here $\phi(\tht)^*$ is obtained from $\phi(\tht)$ by complex conjugation and
reversing the order in the multiplication, i.e.,
\[ (\tht^{i_1}\cdots\tht^{i_k})^*=\bar\tht^{\bar i_k}\cdots\bar\tht^{\bar i_1},
   \quad 1\le i_1<\cdots<i_k\le n,\;0\le k\le n.           \]
The formula bears a formal resemblance with that in
Proposition~\ref{bos-holo}.2 of the bosonic case.
Moreover, the projection from $\HH_0$ to $\HH_J$ is given by the fermionic
counterpart of the Bergman kernel.

\begin{prop}\label{ferm-berg}
The orthogonal projection from $\psi\in\HH_0$ onto $\HH_J$ is 
\[ \psi(\tht,\bar\tht)\longmapsto e^{-\frac{1}{2}g(\tht,\bar\tht)}
   \int_{\Pi V}e^{g(\tht,\bar\chi)-\frac{1}{2}g(\chi,\bar\chi)}
   \,\psi(\chi,\bar\chi)\,\tilde\eps(\chi). \]
\end{prop}

\noindent{\em Proof:}
Suppose the basis is unitary, i.e., $g(e_i,\bar e_{\bar j})=\del_{ij}$.
We write $\tht\bar\chi=\tht^1\bar\chi^{\bar 1}+\cdots+\tht^n\bar\chi^{\bar n}$
for two fermionic vectors $\psi$, $\chi$ in $\Pi V$.
The fermionic measure can be written as 
\[ \tilde\eps(\tht)=d\tht^1 d\bar\tht^{\bar 1}\cdots d\tht^n d\bar\tht^{\bar n}
                   =d\tht d\bar\tht.                      \]
It is easy to check that $\HH_J$ has an unitary basis $\set{\tht^{i_1}\cdots
\tht^{i_k}e^{-\frac{1}{2}\tht\bar\tht}}{0\le k\le n,1\le i_1<\cdots<i_k\le n}$.
So the Bergman kernel that produces the orthoganal projection from $\HH_0$ to
$\HH_J$ is
\[  K(\tht,\bar\chi)=e^{-\frac{1}{2}\tht\bar\tht-\frac{1}{2}\chi\bar\chi}
   \sum_{k=0}^n\sum_{1\le i_1<\cdots<i_k\le n}\tht^{i_1}\cdots\tht^{i_k}
   \bar\chi^{\bar i_k}\cdots\bar\chi^{\bar i_1}
=e^{\tht\bar\chi-\frac{1}{2}\tht\bar\tht-\frac{1}{2}\chi\bar\chi}.  \]
\qed\vspace{10pt}

A fermionic coherent state $c_J^\al$ is of the form
\[   c_J^\al(\tht)=\exp[g(\tht,\bar\al)-\hf g(\tht,\bar\tht)],   \]
where $\al\in\Pi V_J^{1,0}$ is a fermionic parameter.
The above projection can be written as
\[ \psi(\tht,\bar\tht)\longmapsto e^{-\frac{1}{2}g(\tht,\bar\tht)}\int_{\Pi V}
   c_J^\tht(\chi)^*\,\psi(\chi,\bar\chi)\,\tilde\eps(\chi).  \]

Finally, we find the relation to real fermionic coordinates.
Let $\tht_0=\tht+\bar\tht\in\Pi V$.
We have $2\ii g(\tht,\bar\tht)=\varpi_J(\tht_0,\tht_0)$, which we write as
$\varpi_J(\tht_0)$ for short.
Then the fermionic Gaussian factor becomes $e^{\frac{\ii}{4}\varpi_J(\tht_0)}$.

\section{Invariant real, complex and quaternionic structures}\label{C}

\subsection{Representations of real and quaternionic types}

Let $W$ be a finite-dimensional complex vector space.
An operator $C\colon W\to W$ is conjugate linear if $C$ is real linear and
$C(ax)=\bar a(Cx)$ for all $a\in\CO$, $x\in W$.
Such operators are in $\Hom_\CO(W,\bar W)$, where $\bar W$ the complex
vector space which is equal to $W$ as an Abelian group but whose scalar
multiplication is given by $(a,x)\in\CO\times W\mapsto\bar ax\in W$.
A real structure on $W$ is a conjugate-linear operator $R$ on $W$ such that
$R^2=\id_W$.
Such an $R$ determines a real vector space $W_0=W^R\subset W$ fixed by
$R$ and $W\cong(W_0)^\CO$ as complex vector spaces.
A quaternionic structure on $W$ is a conjugate-linear operator $Q$ on $W$
such that $Q^2=-\,\id_W$.
This makes $W$ a quaternionic vector space whith the scalar multiplication
$(a+bj,x)\in\HA\times W\mapsto ax+bQx\in W$ (where $a,b\in\CO$).

Suppose a group $K$ acts on $W$ by a complex representation.
The representation is of real type if there is a $K$-invariant real
structure on $W$.
Such a representation is the complexification of a real representation
of $K$ on $W_0$.
The representation is of quaternionic type if there is a $K$-invariant
quaternionic structure on $W$.
Such a representation is quaternionic-linear with the above scalar
multiplication by $\HA$.
We refer the reader to \cite{BD} for the standard properties of real-
and quaternionic-type representations.
We collect here some more results that will be used in \S\ref{C}.2.

Unless stated otherwise, we assume from now on that $K$ is a finite group
or a compact Lie group.
By averaging over $K$, there is a $K$-invariant Hermitian form
$h\colon W\times\bar W\to\CO$ on $W$.
Our convention of Hermitian forms on $W$ is that they are complex linear
in the first variable but conjugate linear in the second.
(However, we took the opposite convention of physicists for pre-quantum or
quantum Hilbert spaces.)

\begin{lem}\label{rq-repr}
Consider a representation of $K$ on a complex vector space $W$.
Then\\
1. the representation is of real (quaternionic, respectively) type if and
only if there is a non-degenerate $K$-invariant symmetric (skew-symmetric,
respectively) bilinear form on $W$;\\
2. there is a non-zero sub-representation $W'\subset W$ of real (quaternionic, 
respectively) type if and only if there is a non-zero $K$-invariant symmetric
(skew-symmetric, respectively) bilinear form on $W$.
\end{lem}

\noindent{\em Proof:}
Part~1 is well known; see for example Proposition~II.6.4 in \cite{BD} or
the proof of Lemma~\ref{inv-herm} below.
Part~2 follows immediately by taking $W'$ as the orthogonal complement (with
respect to a $K$-invariant Hermitian form) of the kernel of the bilinear form.
\qed\vspace{5pt}

\begin{lem}\label{inv-herm}
Under the same conditions as in Lemma~\ref{rq-repr}, suppose $h$ is a
$K$-invariant Hermitian form on $W$.
If the representation of $K$ on $W$ is of real (quaternionic, respectively)
type, then there is a $K$-invariant real structure $R$ (quaternionic structure
$Q$, respectively)
on $W$ such that $h(Rx,Ry)=h(y,x)$ ($\,h(Qx,Qy)=h(y,x)$, respectively) for all 
$x,y\in W$.
\end{lem}

\noindent{\em Proof:}
Suppose $C_0$ is a $K$-invariant real (quaternionic, respectively)
structure on $W$.
Then $C_0^2=\eps\,\id_W$, where $\eps=\pm1$, respectively.
Consider the complex bilinear form $\beta$ on $W$ given by
$\beta(x,y)=h(x,C_0y)+\eps\,h(y,C_0x)$, where $x,y\in W$.
Then $\beta$ is symmetric (skew-symmetric, respectively) when $\eps=\pm1$.
Moreover, $\beta$ is non-degenerate as $\beta(C_0x,x)=h(x,x)+h(C_0x,C_0x)$
for any $x\in W$.
Following the proof of Proposition~II.6.4 in \cite{BD}, we define an
invertible, conjugate-linear operator $C$ on $W$ by $\beta(x,y)=h(x,Cy)$,
$x,y\in W$.
Then $h(C^2x,y)=\eps\,h(Cy,Cx)=h(x,C^2y)$.
Consequently, $\eps\,C^2$ is self-adjoint and positive definite with respect
to $h$.
The space $W$ decomposes as a direct sum of eigenspaces of $\eps\,C^2$,
each of which is $K$-invariant since $C^2$ is so.
Without loss of generality, assume that $W$ is the eigenspace of $\eps\,C^2$
of a single eigenvalue $\lam>0$.
Since $C^2=\eps\lam\,\id_W$ and $\lam h(x,y)=h(Cy,Cx)$ for all $x,y\in W$,
$\lam^{-1/2}C$ is the desired real (quaternionic, respectively) structure
on $W$ when $\eps=\pm1$.
\qed\vspace{10pt}

The results in Lemma~\ref{inv-herm} can also be explained in matrix language;
we do so when the representation is of real type.
Choosing a real basis of the real subspace $W_0$, the Hermitian form $h$
corresponds to a positive definite Hermitian matrix $H$.
There is a unitary matrix $U$ such that $D={}^T\!UH\bar U$ is a diagonal
matrix of positive entries. 
With the representation of $K$, $U$ can be chosen to commute with $K$.
Let $R=U\,{}^T\!U$.
Since $R\bar R$ is the identity matrix, $R$ defines a real structure on $W$.
The result follows form the identity ${}^T\!RH\bar R=UD{}^T\!\bar U=\bar H$.

\subsection{Complex structures invariant under a representation}

If $V$ is a real vector space and $J$ is a complex structure on $V$, we denote
by $(V,J)$ the complex vector space whose underlying real vector space is $V$
and on which the scalar multiplication by $\ii$ is the action of $J$.
Clearly, $(V,J)\cong V_J^{1,0}$ as complex vector spaces.
Let $K$ be a finite or a compact Lie group.
A presentation of $K$ on $V$ is a complex representation on $(V,J)$
if and only if $J\in\JJ$ is invariant under $K$.

\begin{prop}\label{inv-om-g}
1. Suppose $(V,\om)$ is a symplectic vector space and $J\in\JJ_\om$.
If a representation of $K$ on $(V,J)$ is of real type and preserves $\om$,
then there is a $K$-invariant real structure $R$ on $(V,J)$ such that
$\om(Rx,Ry)=-\om(x,y)$ for all $x,y\in V$.\\
2. Suppose $(V,g)$ is an oriented Euclidean vector space of even dimension
and $J\in\JJ_g$.
If a representation of $K$ on $(V,J)$ is of quaternionic type and preserves
$g$, then there is a $K$-invariant quaternionic structure $Q$ on $(V,J)$
such that $g(Qx,Qy)=g(x,y)$ for all $x,y\in V$.
\end{prop}

\noindent{\em Proof:}
1. Let $h(x,y)=\om(x,Jy)-\ii\om(x,y)$, $x,y\in V$.
Then $h$ is a $K$-invariant Hermitian form, as $h(Jx,y)=-h(x,Jy)=\ii h(x,y)$.
By Lemma~\ref{inv-herm}, there is a $K$-invariant real structure $R$ on $(V,J)$
(i.e., $RJ=-JR$, $R^2=\id_V$) such that $h(Rx,Ry)=h(y,x)$.
This is equivalent to $\om(Rx,Ry)=-\om(x,y)$.\\
2. Let $h(x,y)=g(x,y)-\ii g(Jx,y)$, $x,y\in V$.
Then $h$ is a $K$-invariant Hermitian form, as $h(Jx,y)=-h(x,Jy)=\ii h(x,y)$.
By Lemma~\ref{inv-herm}, there is a $K$-invariant quaternionic structure $Q$
on $(V,J)$ (i.e., $QJ=-JQ$, $Q^2=-\,\id_V$) such that $h(Qx,Qy)=h(y,x)$.
This is equivalent to $g(Qx,Qy)=g(x,y)$.
\qed\vspace{5pt}

Since the group $K$ acts on $\JJ_\om$ ($\JJ_g$, respectively) by isometry,
the fixed-point set $(\JJ_\om)^K$ ($(\JJ_g)^K$, respectively) is a totally
geodesic submanifold.

\begin{prop}\label{K-b-f}
1. Suppose $J_0\in\JJ_\om$ is preserved by a symplectic representation
of $K$ on $(V,\om)$.
Then $T_{J_0}(\JJ_\om)^K\cong(\sym^2(V_{J_0}^{1,0}))^K$.
Moreover, the following statements are equivalent:

(a) $(\JJ_\om)^K$ contains a point other than $J_0$;

(b) $(\sym^2(V_{J_0}^{1,0}))^K\ne\{0\}$;

(c) there is a non-zero complex sub-representation $(V',J_0)$ of $(V,J_0)$
of real type.\\
In this case, $V'$ can be chosen as a symplectic subspace and there
is a $K$-invariant real structure $R$ on $(V',J_0)$ such that
$\om(Rx,Ry)=-\om(x,y)$ for all $x,y\in V'$.\\
2. Suppose $J_0\in\JJ_g$ is preserved by an orthogonal representation of 
$K$ on $(V,g)$.
Then $T_{J_0}(\JJ_g)^K\cong(\medwedge^2(V_{J_0}^{1,0}))^K$.
Moreover, the following statements are equivalent:

(a) $(\JJ_g)^K$ is not a discrete set;

(b) $(\medwedge^2(V_{J_0}^{1,0}))^K\ne\{0\}$;

(c) there is a non-zero complex sub-representation $(V',J_0)$ of $(V,J_0)$ of
quaternionic type.\\
In this case, there is a $K$-invariant quaternionic structure $Q$ on $(V',J_0)$
such that $g(Qx,Qy)=g(x,y)$ for all $x,y\in V'$.
\end{prop}

\noindent{\em Proof:}
1. The result on the tangent space follows from \S\ref{A}.2.
Since $\JJ_\om$ is also $K$-equivariantly diffeomorphic to
$\sym^2(V_{J_0}^{1,0})$, the equivalence of (a) and (b) is clear.
The equivalence with (c) is a consequence of Lemma~\ref{rq-repr}.2 and
the rest follows from Proposition~\ref{inv-om-g}.1.\\
2. Although $(\JJ_g)^K$ is not $(\medwedge^2(V_{J_0}^{1,0}))^K$ globally,
the latter is isomorphic to the tangent space $T_{J_0}(\JJ_g)^K$, which
is zero if and only if $(\JJ_g)^K$ is a discrete set.
The rest of the proof is similar to that of part~1.
\qed\vspace{5pt}

By compactness, the set $(\JJ_g)^K$ is finite if it is discrete.
When $K$ is a compact torus group, the number of elements in $(\JJ_g)^K$
is equal to the Euler characteristic $\chi(\JJ_g)=2^{n-1}$.
(This is the quotient of the order of the Weyl group of $\SO(2n)$ by that
of $\U(n)$.)
For example, $\JJ_g\cong S^2$ when $n=2$.
If $K=S^1$ acts on $\CO^2$ with weights $1$ and $-1$, then the fixed-point
set in $\JJ_g$ is $(\JJ_g)^K=\{\pm J_0\}$, whose cardinal is $2=\chi(S^2)$.

\vspace{20pt}

\noindent {\bf Acknowledgments.}
Part of the work has been reported at various conferences, in 
Oberwolfach (2004), Sendai (2005), Tianjin (2005), Hong Kong (2006).
The author thanks W.~D.~Kirwin, N.~Mok and R.~Sjamaar for helpful
discussions.
The work is supported in part by CERG HKU705407P,  Fondation Sciences
Math\'ematiques de Paris and Universit\'e Paris--Diderot (Paris~7).

\end{document}